\magnification=1200
\font\bsc = cmcsc10 at 16pt
\font\msc = cmcsc10
\font\eightrm =cmr8
\font\eightsl=cmsl8

\newcount\minutes
\newcount\scratch

\def\timestamp{%
\scratch=\time
\divide\scratch by 60
\edef\hours{\the\scratch}
\multiply\scratch by 60
\minutes=\time
\advance\minutes by -\scratch
\the \month/\the\day$\,$---$\,$\hours:\null
\ifnum\minutes< 10 0\fi
\the\minutes}

\centerline {\bsc Tight immersions}
\centerline{\bsc and local differential geometry}
\bigskip
\centerline{\sl Ross Niebergall\footnote{}{\eightrm The work of the
first author was partially
supported by an NSERC grant.} and Gudlaugur Thorbergsson}\bigskip
{\centerline{(December 31, 1996)}

\edef\recoveratcodezqrz{\catcode`\noexpand\@=\the\catcode`\@}

\catcode`\@=11

%

\def\m@ltiplybyscale#1,#2,#3.{%
  \multiply #1 by #2\relax
  \divide #1 by #3\relax
}

\newif\ifAMSTEX@loaded
\newif\ifAMSPPT@loaded

\def\loc@lm@g#1#2#3{%
  %
  \expandafter\ifx\csname amsppt.sty\endcsname\relax
       \AMSPPT@loadedfalse\else\AMSPPT@loadedtrue\fi
  \ifx\amstex@loaded\relax\AMSTEX@loadedtrue\else\AMSTEX@loadedfalse\fi
  %
  %
  \edef\c@ntfive{\the\count5}%
  \count255=#1  
  %
  %
  %
  \ifx\MathTimefontnames\relax
    \message{<Assuming MathTime style font definitions>}%
    %
    \edef\c@untfive{\the\count5}%
    \edef\c@untsix{\the\count6}%
    \edef\c@untseven{\the\count7}%
    \edef\c@unteight{\the\count8}%
    %
    \count5=\count255
    \m@ltiplybyscale\count5,1,2.
    %
    \count6=\count255
    \m@ltiplybyscale\count6,3,5.
    %
    \count7=\count255
    \m@ltiplybyscale\count7,7,10.%
    %
    \count8=\count255
    \m@ltiplybyscale\count8,4,5.
    %
%
    \font\tenrm=\RomanFont scaled \count255
    \font\sevenrm=\RomanFont scaled \count5
    \font\fiverm=\RomanFont scaled \count7
    \textfont0=\tenrm
    \scriptfont0=\sevenrm
    \scriptscriptfont0=\fiverm
%
    \font\teni=MTMI scaled \count255
    \font\seveni=MTMI scaled \count7
    \font\fivei=MTMI scaled \count5
    \ifAMSPPT@loaded
      \font\eighti=MTMI scaled \count8\relax  \skewchar\eighti=45
      \font\sixi=MTMI scaled \count6\relax    \skewchar\sixi=45
    \fi
    \textfont1=\teni
    \scriptfont1=\seveni
    \scriptscriptfont1=\fivei
%
    \font\tensy=MTSY scaled \count255
    \font\sevensy=MTSY scaled \count7
    \font\fivesy=MTSY scaled \count5
    \textfont2=\tensy
    \scriptfont2=\sevensy
    \scriptscriptfont2=\fivesy
%
    \font\tenex=MTEX scaled \count255
    \font\sevenex=MTEX scaled \count7
    \textfont3=\tenex
    \scriptfont3=\sevenex
%
    \font\tenit=\ItalicFont scaled \count255
    \font\sevenit=\ItalicFont scaled \count7
    \textfont\itfam=\tenit
    \scriptfont\itfam=\sevenit
%
    \font\tensl=\SlantFont scaled \count255
    \textfont\slfam=\tensl
%
    \font\tenbf=\BoldfaceFont scaled \count255
    \textfont\bffam=\tenbf
%
    \font\tentt=\TypewriterFont scaled \count255
    \textfont\ttfam=\tentt
    %
    \count5=\c@untfive
    \count6=\c@untsix
    \count7=\c@untseven
    \count8=\c@unteight
  %
  \else    
    \message{<Assuming Computer Modern style font definitions>}%
    \font\tenrm=cmr10 scaled \count255
    \font\sevenrm=cmr7 scaled \count255
    \font\fiverm=cmr5 scaled \count255
    \ifAMSPPT@loaded 
      \font\eightrm=cmr8 scaled \count255
      \font\sixrm=cmr6 scaled \count255
    \fi

    \textfont0=\tenrm
    \scriptfont0=\sevenrm
    \scriptscriptfont0=\fiverm
    \font\teni=cmmi10 scaled \count255
    \font\seveni=cmmi7 scaled \count255
    \font\fivei=cmmi5 scaled \count255
    \ifAMSPPT@loaded
      \font\eighti=cmmi8 scaled \count255
      \font\sixi=cmmi6 scaled \count255
    \fi

    \textfont1=\teni
    \scriptfont1=\seveni
    \scriptscriptfont1=\fivei
    \font\tensy=cmsy10 scaled \count255
    \font\sevensy=cmsy7 scaled \count255
    \font\fivesy=cmsy5 scaled \count255
    \ifAMSPPT@loaded
      \font\eightsy=cmsy8 scaled \count255
      \font\sixsy=cmsy6 scaled \count255
    \fi

    \textfont2=\tensy
    \scriptfont2=\sevensy
    \scriptscriptfont2=\fivesy
    \font\tenex=cmex10 scaled \count255
    \textfont3=\tenex
    \scriptfont3=\tenex
    \scriptscriptfont3=\tenex
  \ifAMSPPT@loaded
      \font\sevenex=cmex7 scaled \count255
      \font\eightex=cmex8 scaled \count255
      \textfont3=\tenex
      \scriptfont3=\sevenex
      \scriptscriptfont3=\sevenex
    \fi

    \font\tenit=cmti10 scaled \count255
    \ifAMSPPT@loaded
      \font\sevenit=cmti7 scaled \count255
      \scriptfont4=\sevenit
      \scriptscriptfont4=\sevenit
      \font\eightit=cmti8 scaled \count255
    \else
      \ifAMSTEX@loaded
        \scriptfont4=\tenit
        \scriptscriptfont4=\tenit
      \fi
    \fi

    \textfont\itfam=\tenit
    \font\tensl=cmsl10 scaled \count255
    \ifAMSTEX@loaded
      \scriptfont5=\tensl
      \scriptscriptfont5=\tensl
    \fi
    \ifAMSPPT@loaded
      \font\eightsl=cmsl8 scaled \count255
    \fi

    \textfont\slfam=\tensl
    \font\tenbf=cmbx10 scaled \count255
    \ifAMSPPT@loaded
      \font\eightbf=cmbx8 scaled \count255
      \font\sixbf=cmbx6 scaled \count255
    \fi

    \textfont\bffam=\tenbf
    \scriptfont\bffam=\sevenbf
    \scriptscriptfont\bffam=\fivebf
    \font\tentt=cmtt10 scaled \count255
    \ifAMSPPT@loaded
      \font\eighttt=cmtt8 scaled \count255
    \fi

    \textfont\ttfam=\tentt
  \fi  
%
%
  \m@ltiplybyscale\baselineskip,#2,#3.\relax
  \m@ltiplybyscale\normalbaselineskip,#2,#3.\relax
  \m@ltiplybyscale\normallineskip,#2,#3.\relax
%
%
  \m@ltiplybyscale\thinmuskip,#2,#3.\relax
  \m@ltiplybyscale\medmuskip,#2,#3.\relax
  \m@ltiplybyscale\thickmuskip,#2,#3.\relax
%
%
  \m@ltiplybyscale\smallskipamount,#2,#3.\relax
  \m@ltiplybyscale\medskipamount,#2,#3.\relax
  \m@ltiplybyscale\bigskipamount,#2,#3.\relax
  \m@ltiplybyscale\jot,#2,#3.\relax
%
%
  \m@ltiplybyscale\parindent,#2,#3.\relax
  \m@ltiplybyscale\parskip,#2,#3.\relax
%
%
%
  \normalbaselines
  \rm
}

\let\localmag\loc@lm@g

\recoveratcodezqrz

%
%
%
\newdimen\rotdimen
\def\vspec#1{\special{ps:#1}}
\def\rotstart#1{\vspec{gsave currentpoint currentpoint translate
   #1 neg exch neg exch translate}}
\def\rotfinish{\vspec{currentpoint grestore moveto}}
%
%
\def\rotr#1{\rotdimen=\ht#1\advance\rotdimen by\dp#1%
   \hbox to\rotdimen{\hskip\ht#1\vbox to\wd#1{\rotstart{90 rotate}%
   \box#1\vss}\hss}\rotfinish}
%
%
\def\rotl#1{\rotdimen=\ht#1\advance\rotdimen by\dp#1%
   \hbox to\rotdimen{\vbox to\wd#1{\vskip\wd#1\rotstart{270 rotate}%
   \box#1\vss}\hss}\rotfinish}%
%
%
\def\rotu#1{\rotdimen=\ht#1\advance\rotdimen by\dp#1%
   \hbox to\wd#1{\hskip\wd#1\vbox to\rotdimen{\vskip\rotdimen
   \rotstart{-1 dup scale}\box#1\vss}\hss}\rotfinish}%
%
%
\def\rotf#1{\hbox to\wd#1{\hskip\wd#1\rotstart{-1 1 scale}%
   \box#1\hss}\rotfinish}%


\def\ni {\noindent}

\def\R{{\bf R}}
\def\C{{\bf C}}
\def\H{{\bf H}}
\def\O{{\bf O}}
\def\Z{{\bf Z}}
\def\K{{\bf K}}
\def\E{{\bf E}}
\def\A{{\bf A}}
\def\conj{\bar}
\def\re{{\rm Re}}
\def\im{{\rm Im}}
\def\o{\omega}
\def\mod{{\rm \ mod }}
\def\p{\partial}
\def\w{\wedge}

   \let\SO=\SOn

  \let\SU=\SUn

  \let\Sp=\Spn

\def\Gtwo{{\mathop{{{\bf G\/}}\kern-.5pt_{{}_2}}}}
\def\Ffour{{\mathop{{{\bf F\/}}\kern-2.5pt_{{}_4}}}}
\def\Esix{{\mathop{{{\bf E\/}}\kern-.5pt_{{}_6}}}}
\def\Eseven{{\mathop{{{\bf E\/}}\kern-.5pt_{{}_7}}}}
\def\Eeight{{\mathop{{{\bf E\/}}\kern-.5pt_{{}_8}}}}

\bigskip

Historically, the study of tight immersions of manifolds had its origins
in the study of immersions with minimal total absolute curvature.
The most significant result in the study of minimal total absolute
curvature immersions is the theorem of Chern and Lashof, which
completely characterizes minimal total absolute curvature immersions
(and tight immersions) of  spheres into a Euclidean space.  An essential
ingredient in this characterization was a reformulation of the problem in
terms of the Morse theory of linear height functions and the topological
characteristics of the manifold being immersed.  In this sense, the
theorem of Chern and Lashof characterizes tight immersions of the
topologically simplest compact manifold.  It is very natural to try and
characterize tight immersions of other manifolds
with topological restrictions.

The most natural candidates are highly connected manifolds; i.e.,
$2k$-dimensional manifolds that are $(k-1)$-connected, but not
$k$-connected.  In the years since the theorem of Chern and Lashof was
proved, tight immersions of these manifolds have been studied
extensively by Kuiper and others. The purpose of the present paper is to
finish a program begun by Kuiper and complete the proof of 
the following theorem: \medskip

\ni{\bf Theorem.} {\sl Let $M$ be a  $2k$-dimensional compact
differentiable manifold that is $(k-1)$-connected, but not
$k$-connected. If $f:M\to \E^N$ is a substantial and tight
immersion into a Euclidean space,
then $N\leq 3k+2$. If $N=3k+2$ then
$k=1,2,4$ or $8$ and $f(M)$ is up to a real  projective
transformation the standard embedding of a projective
plane $\K P^2$, where $\K=\R$, $\C$, $\H$ or $\O$.}\medskip

In the introduction we will give a background to the study of tight immersions
and review what was previously known about tight immersions of 
highly connected manifolds.
\bigskip
\ni {\bf \S0. Introduction}\medskip
\medskip

An immersion $f:M\to \E^N$ is called {\sl tight}  if and only if for almost all
unit vectors $\xi$ the linear height functions $h_\xi:M\rightarrow \R$
given by $h_\xi(p)=\left <\xi,f(p)\right >$ are perfect Morse functions
with respect to $\Z_2$ ({\sl perfect} meaning that the Morse
inequalities are equalities). An immersion is said to be {\sl
substantial} if its image does not lie in any hyperplane of $\E^N$.
The connection between tight immersions and immersions into Euclidean space
with minimal total absolute curvature is explained in
[CL1], [CL2] and also [K4], [CR].

The tightness of an immersion is invariant under real projective
transformations. In fact, only the underlying affine space $\A^N$ is
needed in the definition of a tight immersion, and not the space
$\E^N$ endowed with a Euclidean metric.  Let $H$ be a hyperplane in
the real projective space $\R P^N$ and let $\A^N=\R P^N \setminus H$
be its complement.  Suppose that $f:M\rightarrow \A^N$ is a tight
immersion and $\Gamma:\R P^N\rightarrow \R P^N$ is a projective
transformation such that $\Gamma\circ f(M)\subset \A^N$.  Then
$\Gamma\circ f:M\rightarrow \A^N$ is also tight.

The projective plane $\K P^2$, where $\K$ denotes the real numbers
$\R$, the complex numbers $\C$, the quaternions $\H$ or the octonions
(Cayley numbers) $\O$, is $2k$-dimensional and $(k-1)$-connected, but
not $k$-connected, where $k=\hbox{dim}_\R \K$. Notice that $k=1$, $2$,
$4$ or $8$ for $\K=\R$, $\C$, $\H$ or $\O$ respectively.  The
projective plane $\K P^2$ can be identified with the set of $3\times
3$ matrices $A$ over $\K$ satisfying $A=\bar A^t$, $A^2=A$ and
$\hbox{trace}(A)=1$. We have thus embedded $\K P^2$ as a
differentiable submanifold of the linear space of $3\times 3$ matrices
over $\K$. It turns out that the image is substantial in an affine
subspace with real dimension $3k+2$. This embedding into $\A^{3k+2}$
or the real projective space $\R P^{3k+2}$ that we get by adding a
hyperplane at infinity is called the {\sl standard embedding of $\K
P^2$.} If $\K=\R$ we get the well known real Veronese in $\R P^5$. See
[K4] for a thorough discussion of these examples.

Another way to introduce the standard embeddings of the projective
planes is to look at the isotropy representations of the rank two
symmetric spaces of either compact or non-compact type with Weyl group of
type $A_2$. In the compact type case these spaces are $\SU(3)/\SO(3)$,
$\SU(3)$, $\SU(6)/\Sp(3)$ and ${\bf E}_6/{\bf F}_4$. The exceptional
orbits of the isotropy representation of these symmetric spaces are
diffeomorphic to the projective planes $\K P^2$ where $\K=\R$, $\C$,
$\H$ and $\O$ respectively. In each case these orbits are up to a
homothety congruent to the standard embeddings of the projective planes.
The principal orbits of the isotropy representation of these symmetric
spaces are diffeomorphic to the flag manifolds of the corresponding
projective planes. Bott and Samelson proved in [BS], using a different
terminology, that an orbit of the isotropy representation of a symmetric
space is tight. This was reproved by Kobayashi and Takeuchi in [KT]. The
special case of the standard embeddings of projective planes was
apparently first explicitly considered by Tai in [T], who gave a direct
proof of their tightness.

Tight immersions of $2k$-dimensional $(k-1)$-connected manifolds were
studied extensively by Kuiper.  In fact, our theorem above was proved
by Kuiper in the case $k=1$ in [K2] in 1962, although in this case $M$
is just a closed surface, as the topological conditions simply mean
the manifold is connected. In 1970 in [K3] he proved that $k$ must be
equal to $1$, $2$, $4$ or $8$ if the codimension is at least three,
and that $k+2$ is the highest codimension in which a manifold
satisfying the conditions of the theorem admits a tight substantial
immersion, see also [K1]. In [K4], Kuiper continues studying
immersions of these manifolds, although here he makes the stronger
assumption that $M$ admits a Morse function with exactly three
critical points, or in his terminology, that $M$ is \lq like a
projective plane\rq.  Under this additional assumption, he is able to
prove our above theorem for $k=2$.  He also proves that when $k=4$ or $8$
the immersed submanifold is a smooth algebraic subvariety that carries
a differentiable incidence structure satisfying the axioms of a
projective plane. These investigations of Kuiper were continued in
[Th].

We use local projective differential geometry via the method of moving
frames to prove the theorem. Our approach follows the paper [GH] of
Griffiths and Harris closely. They prove a local version of a theorem
of Severi [S] saying that a piece of a substantial complex
two-dimensional surface in $\C P^5$ with degenerate secant variety and
nondegenerate tangential variety is projectively equivalent to a piece
of the Veronese. Important in their proof is that a refinement
$\widehat{III}$ of the third fundamental form vanishes. This is
somewhat stronger than the property of tight immersions proved by
Little and Pohl in [LP] that a hyperplane supporting the image of the
immersion to the second order also supports it to the third order. We
prove that $\widehat{III}$ vanishes for a tight immersion as in the
theorem if the codimension is $k+2$.  The proof uses a normal form for
the second fundamental form of such tight immersions that is due to
Kuiper, see [K3] and [K1], together with a result from [Th] on
so-called top sets.  Once this is done we can then work strictly in
the framework of local projective differential geometry; there is no
further use for the assumption that the immersion is tight.  In fact
we characterize a piece of a substantial $2k$-dimensional submanifold
of $\R P^{3k+2}$ with vanishing $\widehat{III}$ and second fundamental
form satisfying Kuiper's normal form up to projective equivalence as a
piece of a standard embedding of one of the projective planes. The
topological properties of the underlying manifold only enter in the
proofs of Kuiper's normal form and the vanishing of $\widehat{III}$.
The strategy of our proof is similar to the approach of Little and
Pohl in [LP] or even more to that of Sasaki in [Sa] which was also
inspired by [GH]. The
essential difference between our situation and the one in the papers
[GH], [LP], and [Sa] is that the second fundamental form in their case
is the simplest possible, but in our case it has an intricate
structure that makes the proof much more complicated if $k>1$. The
proof of the theorem in [GH] that we use as a model does not use
anything about the complex numbers and can therefore be used to prove
the corresponding statement about real surfaces in $\R P^5$, but it
does use that the variety is two dimensional, i.e.  $k=1$. The essential part
of our proof will consist in giving a partial generalization of the
real case of the theorem of Griffiths and Harris to the cases $k=2$,
$4$ and $8$.

In Section 1 we review some basic material on local projective differential
geometry from [GH] to fix the 
notation and make the paper more readable. In Section
2 we introduce Kuiper's normal form and prove the vanishing of $\widehat{III}$.
In Section 3 we prove the local theorem which immediately implies the main
result of the paper.
\bigskip

\ni {\bf \S1. Local Differential Geometry}\hfill \medskip

\def\P{{\bf P}}
\def\om{\omega}
\def\we{\wedge}

 In this section we will present the necessary details concerning  the
local differential geometry of submanifolds in projective space
and moving frames. Our
treatment will follow very closely that of Griffiths and Harris [GH].
Although our discussion will concern differentiable submanifolds of
 real projective space, and theirs
concerned complex varieties of complex projective space, the details will be
very  similar.

Throughout this section, we will use $\P^N$ to denote  the real
projective space $\R P^N$.
A {\sl frame}  $\{ A_0, A_1,\ldots ,A_N\}$ for
$\P^N$ is a basis $A_0,A_1,\ldots,A_N$ for $\R^{N+1}$.  The set of all such
frames form a
manifold ${\cal F}(\P^{N})$ that may be identified with ${\bf
GL}_{N+1}$.  If we take $0\leq i,j,k\leq N$, then the structure equations
for a moving frame are
$$\eqalign{dA_i &=\sum_j\bar\om_{ij} A_j\cr
        d\bar\om_{ij}& =\sum_k \bar\om_{ik}\we\bar\om_{kj},}\eqno(1.1)$$
where the $\bar\om_{ij}$ are $1$-forms on ${\cal F}(\P^N)$.
There is a fibering  $$\pi:{\cal F}(\P^N)\longrightarrow \P^N$$ defined
by
$$\pi(\{A_0,A_1,\ldots,A_N\})=[A_0].$$
For any point $p\in \P^N$ the
fibre $\pi^{-1}(p)$ consists of all frames $\{A_0,A_1,\ldots,A_N\}$
whose first vectors are in the equivalence class of $p$, i.e., $p=[A_0]$.  If
we set
$$\bar\om_i=\bar\om_{0i},$$ then the forms
$$\bar\om_1,\;\bar\om_2,\;\ldots,\;\bar\om_N$$ are
horizontal for this fibering meaning that they vanish on the fibres
$\pi^{-1}(p)$.

Let us fix a frame $\{A_0,A_1,\ldots,A_N\}$. Then the $1$-forms
$\bar\om_1,\;\dots\;\bar\om_N$
at $\{A_0,A_1,\ldots,A_N\}$ are lifts under $\pi$
of uniquely defined $1$-forms
$\hat\om_1,\;\dots,\;\hat\om_N$ on
$\P^N$ at $p$ where $p=[A_0]$. In fact, let $U$ be a neighborhood of
$p$ in $\P^N$ and let
$f:U\to {\cal F}(\P^N)$ be a local section with
$f(p)=\{A_0,A_1,\ldots,A_N\}$. Then
$\hat\om_i=f^*\bar\om_i$ and $\hat\om_i$ is unique since
$\pi^*$ is injective. One can show that
$\hat\om_1,\;\dots,\;\hat\om_N$ is a basis of
$T_p^*\P$. Let $v_1,\;\dots\,\;v_N$ denote the dual basis of $T_p\P$.
Then  $v_i$ is tangent to the chord $\overline{A_0 A_i}$ between $[A_0]$
and $[A_i]$ as can be seen from the equation
$$d A_0\equiv\sum_i\bar\om_i A_i \ \ \ \ {\rm mod}\ A_0.$$
Similarly there are unique $1$-forms $\hat\om_{ij}$ that lift to
$\bar\om_{ij}$.

Suppose we are given a submanifold $M^n\subset \P^N$.  Let $\bar M\subset
\R^{N+1}\setminus\{0\}$ denote the preimage of $M$ under the projection
$\R^{N+1}\setminus\{0\}\to
\P^N$. Let $A_0\in \bar M$ be a point lying over $p$. We will view
$T_{A_0}\bar M$ as a linear
subspace of $\R^{N+1}$. Notice that $A_0\in T_{A_0}\bar M$. We then have
$$T_p M \cong T_{A_0}\bar M/\R\cdot A_0$$
where $A_0\in \R^{N+1} -\{0\}$ is any point lying over $p$ and $T_pM$ is
the linear tangent
space of $M$, i.e., $T_pM\subset T_p\P^N$. We can also define the {\sl
projective normal space of
M at
$p\in M$} by
$$N_p M=\R^{N+1}/T_{A_0}\bar M.$$

Associated to $M\subset \P^N$ is the submanifold
${\cal F}(M)\subset {\cal F}(\P^N)$ of {\sl Darboux frames}
$$\{A_0;A_1,A_2,\ldots , A_n;A_{n+1},\ldots, A_N\}$$
defined by the conditions that $A_0$ lies over $p\in M$ and
$\{A_0,A_1,\ldots,A_n\}$ spans $T_{A_0}\bar M$. It then follows that
$\{A_{n+1},\ldots,A_N\}$ mod $T_{A_0}\bar M$ is a basis of $N_p M$.

We shall use the following ranges of indices
$$1\leq \alpha,\beta,\gamma\leq n,\ \ \ \ \ \ \ n+1\leq \mu,\nu\leq N.$$
The $1$-forms $\bar\om_\alpha$ at a Darboux frame $\{A_0;A_1,A_2,\ldots ,
A_n;A_{n+1},\ldots,
A_N\}$, $p=[A_0]$, are the lifts under $\pi^*$ of uniquely defined
$1$-forms $\om_\alpha$ on
$T_{p} M$. Notice that $\om_\alpha$ is the restriction of the $1$-form
$\hat\om_\alpha$
on $T_p\P^N$ to
$T_pM$. More generally, we will denote the restriction of $\hat\om_{ij}$ to
$T_pM$ by $\om_{ij}$.
It is clear from (1.1) that the image of $dA_0$ lies in $T_{A_0}\bar M$
if and only if
$\{\hat\om_\alpha\}$ restricts to 
a basis for $T_p^* M$ and $\hat\om_\mu|T_pM=0$.

It follows from (1.1) that
$$0=d\om_\mu=\sum_\alpha \om_\alpha \we \om_{\alpha \mu}.$$
By the Cartan Lemma
it follows that there exist real valued functions
$q_{\alpha\beta\mu}=q_{\beta\alpha\mu}$ such that
$$\om_{\alpha \mu}=\sum_{\beta} q_{\alpha\beta\mu}\om_\beta.\eqno(1.2)$$
We  define
$$Q_\mu=\sum_{\alpha\beta}q_{\alpha\beta\mu}\om_\alpha\om_\beta,
        \ \ \ \ \ {\rm for}\ \ \mu=n+1,\ldots,N.\eqno(1.3)$$
Then the {\sl projective second fundamental form} is the map
$$II: T_pM\times T_pM \rightarrow N_pM$$
given in coordinates by
$$II(v,w)=\sum_{\alpha\beta\mu}q_{\alpha\beta\mu}
                \om_\alpha(v)\om_\beta(w) A_{\mu}\ \mod\ \ T_{A_0} \bar M$$
where $v=\sum_{\alpha}\om_\alpha(v)v_\alpha$ for $v\in T_p M$,
and $v_1,\ldots,v_n\in T_p M$ is a basis dual to
the basis $\om_1,\ldots,\om_n$ of $T_p^* M$.

By the equation
$$dA_0\equiv\sum_{\alpha}\om_{\alpha}A_\alpha\ \ \ {\rm mod}\ A_0$$
we may write
$${{dA_0}\over{dv_\alpha}}\equiv A_\alpha\ \ \ {\rm mod}\ A_0.$$
Similarly we can reformulate the data contained in
the second fundamental form as follows. First one calculates that
$$d^2 A_0\equiv \sum_{\alpha,\beta,\mu}
        q_{\alpha\beta\mu}\om_{\alpha}\om_{\beta} A_\mu\ \ \
        {\rm mod}\ \ T_{A_0}\bar M.$$
Hence we have
$$II(v_\alpha,v_\beta)\equiv{{\p^2 A_0}\over{\p v_\alpha \p v_\beta}}\equiv
        \sum_{\mu}q_{\alpha\beta\mu} A_\mu\ \ \ {\rm mod}\ \ T_{A_0}\bar M.
        \eqno(1.4)$$

The {\sl second osculating space of $M$ at $p$},  denoted
$T_p^2 M$, is defined  to be the vector space
$$T_p^2 M={\rm span}\left\{ A_0,
        A_\alpha,{{dA_\alpha}\over{dv_\beta}}\right\}
        _{\alpha,\beta=1,\ldots n} \ \ \ \mod \ \ \R A_0,$$
and the {\sl first normal space of $M$ at $p$} is then
$$N_p^1 M=T_p ^2 M/T_p M.$$
The first normal space is the image of the second fundamental form $II$.
Notice that its
dimension is not necessarily constant. The {\sl third fundamental form} can
then be described as
a generalization of equation (1.4) for the second fundamental form.  We set
$$III(v,w,z)={{\p^3 A_0}\over{\p v\,\p w\,\p z}}\ \ \ {\rm mod}\ T^2_p M.$$
This defines a map
$$III: T_pM\times T_pM\times T_pM\longrightarrow \R^{N+1}/T_p^2 M.$$
The following refinement of the third fundamental form will be very
important in the proof of the main result of this paper:
$$\widehat{III}(v_\alpha,v_\alpha,v_\beta)\equiv{{\p^3A_0}\over{\p v_\alpha^2
\p v_\beta}}\ \ \
        {\rm mod\ \ Span}\{T_{A_0}\bar M, II(v_\alpha,T_p M)\}.$$ \medskip

The proof of our main theorem, given in 
Section 3, relies on our
ability to compute the Maurer-Cartan form matrix for the immersions
in which we are interested.  The reason why this is important is revealed
in the following elementary lemma concerning equivalence of
mappings of a connected
manifold $M$ into a Lie group $G$.  We view the Maurer-Cartan forms on $G$
collectively as a $1$-form $\varphi$ with values in the Lie algebra $\cal G$
of $G$.  Given a mapping $f:M\rightarrow G$, the pullback form
$\varphi_f=f^* \varphi$ determines $f$ up to a left 
transformation as follows:

\medskip
\ni{\bf Lemma.}  {\sl Given a pair of maps $f, \tilde f:M\rightarrow G$
there exists a fixed $g\in G$ such that
$$f=g\cdot \tilde f,$$
if and only if $\varphi_f=\varphi_{\tilde f}$.}
\bigskip
\bigskip

\penalty-100
\ni {\bf\S2. Kuiper's normal form and the vanishing of $\widehat{III}$}\medskip
\penalty200
\bigskip
\penalty-100
\ni{\sl \S 2.1 Kuiper's normal form}
\penalty200
\medskip

In this section we will describe a normal form for the second
fundamental form due to Kuiper [K3].
Let $M$ be a manifold as described in the main theorem of
this  paper, then the normal form exists at a 
point $p\in M$ that is a nondegenerate minimum of a
height function. More precisely, we will assume that $M$
is a  $2k$-dimensional compact differentiable manifold which is
$(k-1)$-connected, but not $k$-connected. We assume furthermore
that $f:M\to \E^{3k+2}$ is a substantial
tight immersion. Notice that it is only in the proof of Kuiper's normal
form and in the proof of the vanishing of $\widehat{III}$ that
the
topological conditions on the manifold $M$ enter.

Let $p \in M$ be a nondegenerate minimum of a height function.
Kuiper's normal form allows us to find a Darboux frame
$\{A_0;A_1,\dots, A_{2k};A_{2k+1},\dots, A_{3k+2}\}$ with very special
coefficients $q_{\alpha\beta\mu}$ of the second fundamental form. When
we apply this in the next section to a submanifold in $\P^{3k+2}$ we
will consider $\E^{3k+2}$ to be the hyperplane in $\R^{3k+3}$ such
that the last coordinate of the points is equal to $1$, i.e.,
$\P^{3k+2}$ is $\E^{3k+2}$ with a plane at infinity. 
We therefore have to add $1$ as last
coordinate to $A_0$ and $0$ as a last coordinate in $A_1,\dots,
A_{3k+2}$.  Notice that Kuiper's normal form will hold in a whole
neighborhood of $p$ since the condition of being a nondegenerate
minimum of a height function is open.

We now give a precise statement of the normal form.

There is a frame
$\{A_0;A_1,\dots, A_{2k}; A_{2k+1},\dots, A_{3k+2}\}$ such that 

(i) $A_0=f(p)$,

(ii) $A_1,\dots, A_{2k}$ span $f_*(T_pM)$, 

(iii) $A_{2k+1},\dots, A_{3k+2}$ span the (Euclidean) 
normal space of $f$ at $p$, 

(iv) if $\xi=w_1A_{2k+1}+w_2A_{2k+2}+s_1A_{2k+3}+\dots +s_kA_{3k+2}$, 
then the
matrix of the shape operator $S_\xi:T_pM\to T_pM$
defined implicitly by
$$\langle S_\xi v,w\rangle = \langle II(v,w),\xi\rangle$$
with respect to the
basis $A_1,\dots,A_{2k}$ is
$$\left(\matrix{w_1{\bf I} & B\cr
               ^tB      & w_2{\bf I}\cr}\right)$$ where $^tB$ is the
transpose of a matrix $B$ which is a scalar multiple of an orthogonal
$k\times k$ matrix. 

It follows from the theorem of Hurwitz [H] on the
dimensions of normed algebras that a $k$ dimensional linear family of
multiples of $k\times k$ orthogonal matrices exist if and only if $k=1$,
$2$, $4$ or $8$. Using normal forms for such linear families ([H], [K3]),
one can choose the Darboux frame
such that $B$  is the $k\times k$ submatrix in the
upper left hand corner of the following matrix

\def\back{\noalign{\vskip-9pt}}
\def\sliny{\ \vrule height 6pt depth 5pt}
\def\lliny{\ \vrule height 16pt depth 5pt}
$$\left(\matrix{
s_1\sliny  & -s_2\sliny & -s_3 & -s_4\sliny & -s_5 & -s_6 & -s_7 & -s_8 \cr
\noalign{\hrule width 13pt}
\back
\hfill s_2  & \hfill s_1\lliny & -s_4 & \hfill s_3\lliny & -s_6 &
    \hfill s_5 & -s_8 & \hfill s_7 \cr
\noalign{\hrule width 43.5pt}
\back
\hfill s_3  & \hfill s_4 & \hfill s_1 & -s_2\lliny & -s_7 &
    \hfill s_8 & \hfill s_5 & -s_6 \cr
\back
\hfill s_4  & -s_3 & \hfill s_2 & \hfill s_1\lliny & \hfill s_8
    & \hfill s_7 & -s_6 & -s_5 \cr
\noalign{\hrule width 101pt}
\noalign{\vskip 3pt}
\hfill s_5  & \hfill s_6 & \hfill s_7 & -s_8 & \hfill s_1
    & -s_2 & -s_3 & \hfill s_4 \cr
\noalign{\vskip 3pt}
\hfill s_6  & -s_5 & -s_8 & -s_7 & \hfill s_2 & \hfill
    s_1 & \hfill s_4 & \hfill s_3 \cr
\noalign{\vskip 3pt}
\hfill s_7  & \hfill s_8 & -s_5 & \hfill s_6 & \hfill s_3
    & -s_4 & \hfill s_1 & -s_2 \cr
\noalign{\vskip 3pt}
\hfill s_8  & -s_7 & \hfill s_6 & \hfill s_5 & -s_4 &
    -s_3 & \hfill s_2 & \hfill s_1 \cr
}\right).$$\bigskip

\medskip

\ni{\bf Remark.} 
Notice that the restriction that $k$ be $1$, $2$, $4$
or $8$ follows from the rather elementary theorem of Hurwitz on the
dimensions of normed algebras.  It is not necessary to use deep
theorems in topology due to Adams as in [K3]. The reason is that we
are in the maximal codimension $k+2$. If the codimension is $\ell+2$,
then the argument in the proof of Kuiper's normal form implies that we
have an $\ell$ dimensional linear family of multiples of $k\times k$
orthogonal matrices. It follows that $\ell$ cannot be bigger than
$k$. If $\ell<k$, then some restriction on $k$ follow, but algebraic
arguments alone do not suffice to prove that $k$ is one of the numbers
$1$, $2$, $4$ or $8$.

\bigskip
\penalty-100
\ni{\sl \S 2.2 The vanishing of $\widehat{III}$}
\penalty200
\medskip

Little and Pohl [LP]  proved that if
$f:M\to \E^N$ is a tight immersion and $H$ a hyperplane that supports
$f(M)$ to the second order at a
point
$p\in M$ that is a nondegenerate maximum of a height function, then $H$
supports $f(M)$ to the
third order at that point. We will need a refinement of this result in
terms of the refined third
fundamental form, defined in Section 1,
for tight immersions as in the main theorem of the paper.
More precisely we will
prove:
\medskip

\ni{\bf Proposition.} {\sl
If $f:M^{2k}\to\E^{3k+2}\subset\P^{3k+2}$
is a substantial tight immersion of a
$(k-1)$-connected, but not $k$-connected manifold and assume that $p\in
M$ is a nondegenerate minimum of a height function.
Then
$$\widehat{III}(v,v,w)=0$$
for all $v,w\in T_pM$.
}\medskip

\ni{\sl Proof.}
Let $v,\ w\in T_pM$.
We will then show that
the following holds: There is  a map $g: U\to M$, where $(0,0)\in U\subset
\R^2$, such that the
partial derivatives of $g$ satisfy
$$g_x(0,0)=v \ \ \ \hbox{ and }\ \ \ g_y(0,0)=w$$
and the partial derivative $G_{xxy}(0,0)$ lies in the space spanned by
$T_pM$ and
$II(v,T_pM)$, where $G=f\circ g$. This implies the claim in the lemma.

It follows from Kuiper's normal form that $II(v,T_pM)$ is a
$(k+1)$-dimensional subspace of the (Euclidean) normal space.
Hence it follows that the span $H$ of $T_pM$ and $II(v,T_pM)$
is a hyperplane. It follows furthermore from the normal form that $II(v,T_pM)$
supports the image $\{II(w,w)\ \vert \ w\in T_pM\}$ of the second fundamental
form in a ray. By [Th, p.111] it follows that

(i) $H$ supports $f(M)$,

(ii) the boundary $Q$ of the convex hull of $f(M)\cap H$ is contained in $f(M)$
and $f(M)\cap H$ spans a $(k+1)$-dimensional affine subspace of $\E^{3k+2}$
(notice
that $f(M)\cap H$ is a so-called top set, i.e., $f(M)\cap H$ is the set
of minima of a height function, see (i)),

(iii) $f$ is injective on $f^{-1}(Q)$,

(iv) $f^{-1}(Q)$ is a non-trivial singular $k$-cycle, and

(v) there is a neighborhood $V$ of $p$ in $M$ such that $V\cap f^{-1}(H)$
is a $k$-dimensional differentiable submanifold with 
$v\in T_p(V\cap f^{-1}(H))$.

We will only use properties (i) and (v) above. The other
properties should make it easier for the reader 
to use the reference [Th]. One can continue this
study and show that the set $Q$ is actually a quadric. This was done by
Kuiper in [K4], but we will not need such detailed information here.

We can now find a map $g:U\to M$ satisfying the condition that
$g(x,0)\in V\cap f^{-1}(H)$ for all
$(x,0)\in U$ and $g_x(0,0)=v$ and $g_y(0,0)=w$. Let
$\xi$ be perpendicular to $H$ and such that the height function $h_\xi:M\to
\R;q\to \langle
\xi,f(q)\rangle$ has a minimum in $p$. Set $k=h_\xi\circ g$. Then $k:U\to
R$ satisfies $k\ge 0$,
$k(x,0)=0$ for all $(x,0)\in U$ and $k_{xx}(0,0)=0$. Clearly
$k_{xy}(0,0)=0$ and $k_{xxx}(0,0)=0$.
We look at the Taylor expansion
$$k(x,y)={1\over{ 2 !}}k_{yy}(0,0)\,y^2+{1\over
{3!}}\big(3k_{xxy}(0,0)\,x^2y+3k_{xyy}(0,0)\,xy^2+k_{yyy}(0,0)\,y^3\big)+
R(x,y).$$
It now follows that $0=k_y(x,0)={1\over 2}k_{xxy}(0,0)\,x^2+R_y(x,0)$. 
This in turns implies that
$R_{xxy}(0,0)=-k_{xxy}(0,0)$. But $R_{xxy}(0,0)=0$ so we have proved that
$k_{xxy}(0,0)=0$. It
follows that $G_{xxy}(0,0)$ lies in the hyperplane $H$ which coincides with
the linear span of $T_pM$ and
$II(v,T_pM)$ as we have already observed. This finishes the
proof of the proposition.
\hfill $\sqcap \!\!\!\!\sqcup$
\bigskip

\ni{\bf \S3. Proof of the Theorem}\medskip

In this section we complete the proof of the theorem described at the
beginning of the paper.  We know from Section 2
that if $f:M^{2k}\rightarrow \E^{3k+2}\subset\P^{3k+2}$ 
is a tight substantial immersion of a 
manifold as described in the main theorem, then $k=1$, $2$, $4$ or $8$, 
and at a points in $M$ that are a nondegenerate minima of some height
function we are provided with Kuiper's normal form 
and the refined third fundamental form $\widehat{III}$ vanishes.
The part of our
theorem that remains unproven is that for each of the cases
$k=1$, $2$, $4$, and $8$ the embedding is one of the standard embeddings of
the projective planes up to a projective transformation.  
The proof is done separately for each of the four values of $k$.

In [GH] Griffiths and Harris prove that if $M$ is a piece of a complex
surface immersed in $\C P^5$ and has degenerate secant variety and
nondegenerate tangential variety, then it must be projectively
equivalent to a piece of the Veronese embedding of the complex
projective plane.  Although their proof is for a complex surface in
complex projective space the calculations do not depend on the field
of complex numbers, and can be used to prove our $k=1$ case of a real
surface in $\R P^5$.  The important point in the proof of Griffiths
and Harris is that a surface with degenerate secant variety and
nondegenerate tangential variety has a vanishing refined third
fundamental form $\widehat{III}$. That, along with a normal form for
the second fundamental form is sufficient to prove their theorem.
Hence a tight surface in five space is projectively equivalent to a
piece of the real Veronese in a neighborhood of a point $p$ that is a
nondegenerate minimum of a height function. It is now easy to show
that all points are nondegenerate minima of height functions and thus
finish the proof for $k=1$, see Section 3.10 below.

Although the proofs of each of the other cases share similarities with
the proof of Griffiths and Harris there are also significant
differences, not in the least of ways by the increase in magnitude and
complexity of the calculations involved.  The primary goal of the
proofs in each of the cases is the same, and that is to make changes
of frame that make the matrix of Maurer-Cartan forms equivalent to the
matrix corresponding to one of the standard embeddings, and then apply
the Lemma at the end of Section 1. 
For values of $k>1$, we will see that the matrix of
Maurer-Cartan forms has additional structure not present in the $k=1$
case.  This results in the need for additional changes of frame, and
caution must be observed with regards to the order in which these
changes of frame are made.  Since the calculations for the three
remaining cases are very similar, we will describe the proof for the
$k=2$ case in detail and make comments at the end of this section
concerning the cases $k=4$ and $k=8$.

We will divide the proof in this case into several steps.  
We first describe the standard embedding of $\C
P^2$ in $\P^8$ and the matrix of Maurer-Cartan forms associated with
that embedding.  Then let $M^4\subset\P^8$ be an arbitrary tightly
immersed substantial $1$-connected submanifold.  In Section 2 we gave
a precise description of Kuiper's normal form for such an immersion
that is valid in a neighborhood of a point $p$ that is a nondegenerate
minimum of a height function.  The normal form gives us detailed
information concerning the Maurer-Cartan forms for this immersion.  By
carefully analyzing this information, making repeated use of
Cartan's Lemma, and making use of the vanishing
of $\widehat{III}$ proven in the Proposition in Section 2, we are able
to determine significant information concerning the matrix of
Maurer-Cartan forms.  We are then able to make several changes of
frame, solving several systems of linear equations, until finally the
matrix of Maurer-Cartan forms for this arbitrary tight submanifold
agrees with that of the standard embedding $\C P^2\subset
\P^8$.$^\dagger$
\footnote{}{{\vbox{\hsize 6.1truein
\baselineskip 9pt\eightrm  $^\dagger$The
calculations necessary to make the frame changes are rather
formidable; even in the {\eightsl k}=2 case there are several systems
of hundreds of linear equations to be solved.  We used a computer
algebra system to assist in solving these equations.}}}  
Then by the
Lemma in Section 1 we know that a neighborhood around $p$ of our tight
substantially immersed manifold $M^4\subset \P^8$ is equivalent to a
piece of the standard embedding of $\C P^3\subset \R^8$ up to a
projective transformation.
                                                        
\bigskip
\penalty-100
\ni{\sl \S 3.1 The standard embedding of $\C P^2$ in $\P^8$}
\penalty200
\medskip

We first describe the standard embedding of $\C P^2$ in $\P^8$.  This 
description is projectively equivalent to the description given previously.
Let $\{X_0,X_1,X_2\}$ be a basis for $\C^3$.
If we write $X_j=C_j+iD_j$, where $i=\sqrt{-1}$ then we can define the map 
$F:\C^3\rightarrow \P^8$ by 
$$F(X_0,X_1,X_2)=(A_0,A_1,A_2,A_3,A_4,A_5,A_6,A_7,A_8),$$
where
$$\eqalign{A_0&=X_0\conj X_0=C_0^2+D_0^2\cr
           A_1&=\re(X_0 \conj X_1) = C_0 C_1+D_0 D_1\cr
           A_2&=\im(X_0\conj X_1) = C_1 D_0 -C_0 D_1\cr
           A_3&=\re(X_0\conj X_2) = C_0 C_2+D_0 D_2\cr
           A_4&=\im(X_0\conj X_2) = C_2 D_0 -C_0 D_2\cr
           A_5&=X_1\conj X_1 = C_1^2 +D_1 ^2\cr
           A_6&=\re(X_1\conj X_2)=C_1 C_2+D_1 D_2\cr
           A_7&=\im(X_1\conj X_2)=C_2 D_1-C_1 D_2\cr
           A_8&= X_2\bar X_2 =C_2^2+D_2^2.}$$
We can then define the Maurer-Cartan forms for this embedding as follows.
Let $\theta_{jk}=\alpha_{jk}+i\beta_{jk}$, then
$$\eqalign{dX_j&=\sum_{k=0}^{3}\theta_{jk}X_k \cr
                &=\sum_{k=0}^{3}(\alpha_{jk}+i\beta_{jk})(C_k+iD_j)\cr
                &=\sum_{k=0}^3 (\alpha_{jk}C_j-\beta_{jk}D_k+
                        i(\beta_{jk}C_k+\alpha_{jk}D_k)).}$$
Since $dX_j=dC_j+idD_j$, then
$$\eqalign{dC_j=\sum(\alpha_{jk}C_k-\beta_{jk}D_k)\cr
           dD_j=\sum (\beta_{jk}C_k+\alpha_{jk}D_k).}$$
Then the Maurer-Cartan form matrix is given by 
$[dA_0,dA_1,\ldots,dA_8]^T$.  This can be calculated explicitly to obtain

{\localmag{900}{9}{10}
$$ 
\left(
\matrix{ 2\,\alpha_{00} & 2\,\alpha_{1} & 2\,\beta_{1} & 2\,\alpha
_{2} & 2\,\beta_{2} & 0 & 0 & 0 & 0 \cr 
\alpha_{10} & \alpha_{11}+
\alpha_{00} & \beta_{11}-\beta_{00} & \alpha_{12} & \beta_{12} & \alpha
_{1} & \alpha_{2} & \beta_{2} & 0 \cr 
-\beta_{10} & \beta_{00}-\beta_{
11} & \alpha_{11}+\alpha_{00} & -\beta_{12} & \alpha_{12} & \beta_{1}
 & \beta_{2} & -\alpha_{2} & 0 \cr 
\alpha_{20} & \alpha_{21} & \beta_{
21} & \alpha_{22}+\alpha_{00} & \beta_{22}-\beta_{00} & 0 & \alpha_{1}
 & -\beta_{1} & \alpha_{2} \cr 
-\beta_{20} & -\beta_{21} & \alpha_{21}
 & \beta_{00}-\beta_{22} & \alpha_{22}+\alpha_{00} & 0 & \beta_{1} & 
\alpha_{1} & \beta_{2} \cr 
0 & 2\,\alpha_{10} & -2\,\beta_{10} & 0 & 0
 & 2\,\alpha_{11} & 2\,\alpha_{12} & 2\,\beta_{12} & 0 \cr 
0 & \alpha
_{20} & -\beta_{20} & \alpha_{10} & -\beta_{10} & \alpha_{21} & \alpha
_{22}+\alpha_{11} & \beta_{22}-\beta_{11} & \alpha_{12} \cr 
0 & -\beta
_{20} & -\alpha_{20} & \beta_{10} & \alpha_{10} & -\beta_{21} & \beta
_{11}-\beta_{22} & \alpha_{22}+\alpha_{11} & \beta_{12} \cr 
0 & 0 & 0
 & 2\,\alpha_{20} & -2\,\beta_{20} & 0 & 2\,\alpha_{21} & -2\,\beta_{
21} & 2\,\alpha_{22} \cr}
\right)\eqno(3.1)$$}
where the $(j,k)$-th entry is the $A_k$-th component of $dA_j$.  Notice
that $0\leq j,k\leq 8$.
\medskip

\ni{\bf Remark.} It will often be necessary for us to refer to specific
locations and rectangles of entries of this matrix.  We will use the
convention that the rows and columns of the matrix are numbered 
from $0$ to $8$, and we will denote the entry in the $j$-th row and $k$-th
column by $(j,k)$.

\bigskip
\penalty-100
\ni{\sl \S 3.2.  The second fundamental form of a tight $1$-connected
                manifold $M^4$ in $\P^8$}
\penalty200
\medskip

Now we will use the information that we have concerning a tight
substantial immersion of a $1$-connected manifold 
$M^4\subset \P^8$ to show that after
several changes of frame, it has the same Maurer-Cartan form
matrix as (3.1).  We will use the formalism established in Section 1.  
Let $\{A_0;A_1,\ldots,A_4;A_5,\ldots,A_8\}$ be any
Darboux  frame for the immersion, with position vector $A_0$. We denote
the matrix of Maurer-Cartan forms for this immersion by
$\Omega=(\o_{jk})$ where we write $\o_{j}=\o_{0j}$. It follows that
$\o_j=0$ for $j=5,\ldots,8$. We let $p$ be a point that is a nondegenerate
minimum of a height function. Then $p$ has a neighborhood $U$ on which the
same property holds. We will assume that our frames are defined on such
a neighborhood.
The initial piece of
information that we have concerning the immersion restricted to $U$
is Kuiper's normal form for 
the second fundamental form, given in Section 2.1. 
In the case of $k=2$, the matrix of this
form is
$$S=\pmatrix{ w_1 & 0 & s_1 & -s_2 \cr
            0 & w_1 & s_2 & s_1 \cr
            s_1 & s_2 & w_2 & 0 \cr
            -s_2 & s_1 & 0 & w_2 \cr
        }.
$$ 
If we choose the Darboux frame to be compatible with this normal form,
then in the notation of Section 1 we have  
$Q_\mu={{1}\over{2}}[\o_1\o_2 \o_3 \o_4]S[\o_1\o_2 \o_3 \o_4]^T$ from 
equation (1.3).
For $j=1,2$, we compute:
$$\eqalign{{\rm for }\  w_1=1,\ w_2=0,\ s_j=0 :\ \ \ 
                Q_5&= {{1}\over{2}}(\o_1^2+\o_2^2),\cr
           {\rm for }\  s_1=1,\ s_2=0,\ w_j=0 :\ \ \ 
                Q_6&= {{1}\over{2}}(\o_1\o_3+\o_3\o_1+\o_2\o_4+\o_4\o_2),\cr
           {\rm for }\  s_2=1,\ s_1=0,\ w_j=0 :\ \ \ 
                Q_7&= {{1}\over{2}}(\o_1\o_4+\o_4\o_1-\o_2\o_3-\o_3\o_2),\cr
           {\rm for }\  w_2=1,\ w_1=0,\ s_j=0 :\ \ \ 
                Q_8&= {{1}\over{2}}(\o_3^2+\o_4^2).}\eqno(3.2)$$
Let $1\leq \alpha,\beta \leq 4$ and $5\leq \mu \leq 8$. 
By (1.2) we know that $\o_{\alpha\mu}=\sum_{\beta}q_{\alpha
\beta\mu}\o_{\beta}$, where the 
$q_{\alpha\beta\mu}$ are given by the
equation  
$Q_\mu=\sum_{\alpha,\beta}q_{\alpha\beta\mu}\o_\alpha \o_\beta$.
These relations allow us to relate several of the forms in our
matrix.

\def\ha{{{1}\over{2}}}
\def\ss{\hskip 0.5truecm}
\def\vs{\noalign{\vskip 5pt}}
$$\matrix{ \o_{15}=\ha \o_1,\ss & \o_{16}=\ha \o_3,\ss &
                \o_{17}=\ha\o_4,\ss & \o_{18}=0,\cr
\vs
            \o_{25}=\ha \o_2,\ss & \o_{26}=\ha \o_4,\ss &
                \o_{27}=-\ha\o_3,\ss & \o_{28}=0,\cr
\vs
           \o_{35}=0,\ss & \o_{36}=\ha \o_1,\ss &
                \o_{37}=-\ha\o_2,\ss & \o_{38}=\ha \o_3,\cr
\vs
           \o_{45}=0,\ss & \o_{46}=\ha \o_2,\ss &
                \o_{47}=\ha\o_1,\ss & \o_{48}=\ha \o_4.\cr}\eqno(3.3)$$

These calculations then show that the matrix $\Omega$ has a form
identical to that of (3.1) (by scaling by a factor of 2) 
along the top row and in the rectangle with corners at $(0,5)$
and $(4,8)$. Hence, the information that we now
have concerning our matrix $\Omega$ results in the following.

$$ \pmatrix{ 
\omega_{00} & \omega_{1} & \omega_{2} & \omega_{3}
 & \omega_{4} & 0 & 0 & 0 & 0 \cr
\omega_{10} & \omega_{11} & \omega_{12} & \omega_{13}
 & \omega_{14} & {{\omega_{1}}\over{2}} & {{\omega_{3}}\over{2}} & 
{{\omega_{4}}\over{2}} & 0 \cr 
\omega_{20} & \omega_{21} & \omega_{22} & \omega_{23}
 & \omega_{24} & {{\omega_{2}}\over{2}} & {{\omega_{4}}\over{2}} & 
-{{\omega_{3}}\over{2}} & 0 \cr 
\omega_{30} & \omega_{31} & \omega_{32} & \omega_{33}
 & \omega_{34} & 0 & {{\omega_{1}}\over{2}} & -{{\omega_{2}}\over{2}} 
& {{\omega_{3}}\over{2}} \cr 
\omega_{40} & \omega_{41} & \omega_{42} & \omega_{43}
 & \omega_{44} & 0 & {{\omega_{2}}\over{2}} & {{\omega_{1}}\over{2}} 
& {{\omega_{4}}\over{2}} \cr 
\omega_{50} & \omega_{51} & \omega_{52} & \omega_{53}
 & \omega_{54} & \omega_{55} & \omega_{56} & \omega_{57} 
& \omega_{58} \cr 
\omega_{60} & \omega_{61} & \omega_{62} & \omega_{63}
 & \omega_{64} & \omega_{65} & \omega_{66} & \omega_{67} 
& \omega_{68} \cr 
\omega_{70} & \omega_{71} & \omega_{72} & \omega_{73}
 & \omega_{74} & \omega_{75} & \omega_{76} & \omega_{77} 
& \omega_{78} \cr 
\omega_{80} & \omega_{81} & \omega_{82} & \omega_{83}
 & \omega_{84} & \omega_{85} & \omega_{86} & \omega_{87} 
& \omega_{88} \cr
}$$

\bigskip
\penalty-100
\ni{\sl \S 3.3.  Setting up the frame}
\penalty200
\medskip

We can now begin the work of trying to exploit the information contained
in equations (3.3).
The forms $\o_1$, $\o_2$, $\o_3$, and $\o_4$
form a basis for the $1$-forms of $M$.  Thus for each
value of $j,k$ in the range $1\leq j \leq 8$ and $0\leq k \leq 8$ and
also for $j=k=0$ there are functions $b_{jk}^i:U\rightarrow \R$ defined
on an open neighborhood $U$ of $p$ in $M$ such that
$$\o_{jk}=b_{jk}^1 \o_1+b_{jk}^2 \o_2+b_{jk}^3 \o_3+b_{jk}^4\o_4.
        \eqno(3.4)$$     

Recall the set of equations (3.3).  By taking the exterior derivative
of each of these equations, 
we arrive at additional relationships
between the forms of the matrix $\Omega$.  For example since
$$2d\o_{15}=\o_{11}\w\o_1+\o_{12}\w\o_2+
        \o_1\w\o_{55}+\o_{3}\w\o_{65}+\o_{4}\w\o_{75}$$
$$d\o_1=\o_{00}\w\o_{1}+\o_{1}\w\o_{11}+\o_{2}\w\o_{21}
                +\o_{3}\w\o_{31}+\o_{4}\w\o_{41},$$
then the first of equations (3.3), $\o_{15}={{1}\over{2}}\o_1$ 
could be written as
$$(\o_{00}-2\o_{11}+\o_{55})\w\o_1-(\o_{21}+\o_{12})\w \o_{2}
        +(\o_{65}-\o_{31})\w\o_3+(\o_{75}-\o_{41})\w\o_4=0.\eqno(3.5)$$
By replacing each of the non-basis forms
$\o_{jk}$ with its expansion in terms of basis forms as in (3.4), we can
rewrite equation (3.5) completely in terms of wedges of basis forms,
$\o_j\w\o_k$.  Now by collecting the coefficients of each of these
terms, and equating them to zero, we get linear equations in terms of
the coefficients $b_{jk}^l$.  In the case of (3.5) this would give 
us the following equations:
$$
\matrix{
{\rm Coefficient\ of}\ \o_1\w\o_2:&\ \ \ \
        -b_{00}^2+2b_{11}^2-b_{55}^2-b_{21}^1-b_{12}^1=0\cr
{\rm Coefficient\ of}\ \o_1\w\o_3:&\ \ \ \
        -b_{00}^3+2b_{11}^3-b_{55}^3+b_{65}^2-b_{31}^2=0\cr
{\rm Coefficient\ of}\ \o_1\w\o_4:&\ \ \ \
        -b_{00}^4+2b_{11}^4-b_{55}^4+b_{75}^1-b_{41}^1=0\cr
{\rm Coefficient\ of}\ \o_2\w\o_3:&\ \ \ \ 
        b_{21}^3+b_{12}^3+b_{65}^2+b_{31}^2=0\cr
{\rm Coefficient\ of}\ \o_2\w\o_4:&\ \ \ \ 
        b_{21}^4+b_{12}^4+b_{75}^2-b_{41}^2=0\cr
{\rm Coefficient\ of}\ \o_3\w\o_4:&\ \ \ \ 
        -b_{65}^4+b_{31}^4+b_{75}^3-b_{41}^3=0.}\eqno(3.6)$$
Taking the exterior derivative of each of the sixteen equations 
in (3.3), and expanding the result in terms of the basis forms in this way
results in a set of 96  linear equations, some of which may be zero.  
This set does not contain enough equations to determine all the 
coefficients $b_{jk}^i$ that appear in the equations,
but it is possible to solve for some of the coefficients 
arising in terms of other of these coefficients.  For example, we
could solve the first equation in (3.6) for say $b_{00}^2$, to obtain
$b_{00}^2=2b_{11}^2-b_{55}^2-b_{21}^1-b_{12}^1$.
This gives significant information concerning forms in the 
matrix $\Omega$.

It should be noted that the information 
obtained in this way from the equations
(3.3) is independent of the particular frame chosen, in the sense that
any frame that results in Kuiper's normal form will suffice.
It should also be noted that the information  
that we obtain from equation (3.3) in this way is not sufficient to 
characterize any of the forms in the matrix $\Omega$ completely. Still, it 
does provide relationships between different forms of the matrix.

\bigskip
\penalty-100
\ni{\sl \S 3.4. Consequences of the vanishing of $\widehat{III}$}
\penalty200
\medskip

\def\p{\partial}
In the Proposition in Section 2 we saw that 
$\widehat{III}(v,v,w)=0$ for all $v,w\in T_{A_0} \bar{M}$.   Therefore
$$\widehat{III}(v_1,v_1,v_1)
        \equiv{{\p^3 A_0}\over{\p v_1^3}}\equiv0 \mod\ \{A_0,\ldots,A_7\}.$$
But using (3.2) we can also calculate that
$${{\p ^3 A_0}\over{\p v_1^3}}=\ha {{\p A_5}\over{\p v_1}} \mod
\ \{A_0,\ldots,A_7\}.$$
Therefore $\o_{58}(v_1)=0$.  Performing a similar calculation using
the fact that 
$$\widehat{III}(v_1,v_1,v_2)
\equiv{{\p^3 A_0}\over{\p v_1^2 \p v_2}}\equiv 0\ \mod\ \{A_0,\ldots,A_7\}$$
we can calculate that $\o_{58}(v_2)=0$.  By performing similar
calculations for $v_3$ and $v_4$ we see that $\o_{58}=0$.  Nearly
identical equations also show that $\o_{85}=0$.

Next, let $w=v_1+v_3$.  Then
$${{\p A_0}\over{\p w}}={{\p A_0}\over{\p v_1}}+{{\p A_0}\over{\p
        v_3}}=A_1+A_3   \mod\  A_0,$$
and
$$\eqalign{{{\p^2 A_0}\over{\p w^2}}&={{\p A_1}\over{\p v_1}}+{{\p
                A_1}\over{\p v_3}}+{{\p A_3}\over{\p v_1}}+{{\p
                A_3}\over{\p v_3}} \mod\ A_0 \cr
                &={{1}\over{2}}(A_5+2 A_6+A_8) \mod\ A_0, A_1,\ldots, A_4.}$$
Then, we compute that
$$\eqalign{0&= {{\p^3 A_0}\over{\p v_i \p w^2}}={{1}\over{2}}
        \left( {{\p A_5}\over{\p v_i}} +2{{\p A_6}\over{\p v_i}} 
        +{{\p A_8}\over{\p v_i}}\right )\cr     
        &={{1}\over{2}} \big( \o_{55}(v_i)A_5
                +\o_{56}(v_i)A_6+\o_{57}(v_i)A_7 +\o_{58}(v_i)A_8\cr
        &\ \ \ \ \ \ +2\o_{65}(v_i)A_5 +2\o_{66}(v_i)A_6+2\o_{67}(v_i)A_7
                +2\o_{68}A_8\cr
        &\ \ \ \ \ \ +\o_{85}(v_i)A_5+\o_{86}(v_i)A_6 +\o_{87}(v_i)A_7
                +\o_{88}(v_i)A_8 \big)\cr 
                &\mod\ A_0,A_1,\ldots,A_4, 
                II(w,T_{A_0}\bar{M}),
        }\eqno(3.7)$$
where $II(w,T_{A_0}\bar{M})$ is given by the set of vectors
$$II(v_1+v_3,\alpha v_1+\beta v_2+\gamma v_3+\delta v_4)=
                {{1}\over{2}}\big( \alpha A_5 +(\alpha +\gamma) A_6
                                +(\delta - \beta)A_7 +\gamma A_8\big).
                \eqno(3.8)$$
Equation (3.8) can then be used to reduce (3.7) to see that
$$\o_{55}+2\o_{65}+2\o_{68}+\o_{88}=\o_{56}+2\o_{66}+\o_{86}.\eqno(3.9)$$

By varying the choice of $w$ we can obtain other equations like (3.9).
The choices of $w$ that give distinct equations are $w=v_1+v_3$,
$w=v_1+v_4$, and $w=v_2+v_3$.  The complete set of equations that we
obtain from the vanishing of $\widehat{III}$ are 
$$\o_{58}=0\ \ \ \ \ \ \o_{85}=0\eqno(3.10)$$
$$\eqalign{\o_{55}+2\o_{65}+2\o_{68}+\o_{88}&=\o_{56}+2\o_{66}+\o_{86}\cr
        \o_{55}+2\o_{75}+2\o_{78}+\o_{88}&=\o_{57}+2\o_{77}+\o_{87}\cr
        \o_{55}-2\o_{75}-2\o_{78}+\o_{88}&=-\o_{57}+2\o_{77}-\o_{87}
                }\eqno(3.11)$$ 

While (3.10) above gives us immediate information concerning the matrix
$\Omega$, the equations (3.11) contain more subtle information. By
expanding the forms in (3.11) in terms of the basis forms using (3.4),
and then equating the coefficients of each of the basis forms we 
obtain a set of twelve linear equations involving the coefficients
$b_{ij}^k$. Once again, we can solve equations in this collection as in
Section 3.3.

\bigskip
\penalty-100
\ni{\sl \S 3.5.  First and Second Changes of Frame}
\penalty200
\medskip

Since we are trying to show that our Maurer-Cartan forms matrix $\Omega$
can  be put into the same form as (3.1) it is useful to closely examine
(3.1) and identify features that  our matrix $\Omega$ must have. One of
the first things that should be  noticed about the matrix (3.1) is that
certain square minors in this matrix have an
{\sl anti-symmetric} structure.  By that
we mean that the minor is a linear combination of the identity and
a skew-symmetric matrix.
In particular, the square
minors with corners  at  (1,3) and (2,4), and also at (3,1)
and (4,2)  both have this anti-symmetric
structure.   (There are other location where this structure occurs, but
for now we will only address these two locations.)   The
linear systems of equations in Sections 3.3 and 3.4  
do not imply
that our matrix $\Omega$ has this anti-symmetric structure.  It
does however give us sufficient information to 
make a change of frame that achieves this.

As noted, the frames $\{A_0,\ldots,A_8\}$ are not uniquely determined by
(3.3), (3.10) and (3.11). A change of frame of the form 
$$\tilde A_j=A_j+a_j A_0,\ \  {\rm for }\ \ j=1,2,3,4\eqno(3.12)$$ 
for $a_j:U\rightarrow \R$ defined on an open neighborhood $U$ of $p$ 
preserves the  restrictions on the Darboux frames
with respect to Kuiper's normal form. The collection of equations (3.3),
(3.10) and (3.11) remain unaltered by such  a change of frame. The goal
is then to make a change of frame of the form (3.12) which  will result
in a new Darboux frame  $\{A_0,\tilde A_1,\ldots,\tilde
A_4,A_5,\ldots,A_8\}$, such that the matrix of Maurer-Cartan forms
associated with this frame has the anti-symmetric form in the minors
with corners at  $(1,3)$ and $(2,4)$ as well as $(3,1)$ and $(4,2)$.
Although a change of frame of this type will alter the forms in the 
matrix $\Omega$, it will preserve the relations in $\Omega$ provided
by Kuiper's normal form.
Within the rectangle with corners at (1,1) and (4,4)  the
forms will be transformed as \def\ot{\tilde \omega}
$$\ot_{jk}=\omega_{jk}+a_j\omega_{k}.$$ Hence, to obtain the
anti-symmetric form in the minors cornered at locations $(1,3)$ and
$(2,4)$, and also at $(3,1)$ and $(4,2)$  we need to have
$$\eqalign{\ot_{13}&=\ot_{24}\cr
           \ot_{14}&=-\ot_{23}\cr
           \ot_{31}&=\ot_{42}\cr
           \ot_{41}&=-\ot_{32}.}\eqno(3.13)$$
In fact, it follows from the linear equations determined in Sections 
3.3 and 3.4, that it
is possible to make a change of frame which satisfies (3.13).
This change of frame is given by (3.12) with the coefficients $a_j$
given by  
$$\eqalign{a_1&=b_{00}^1+2b_{77}^1-b_{88}^1-2b_{22}^1-b_{77}^3+b_{66}^3\cr
          a_2&=b_{41}^1-2b_{77}^4-b_{00}^4+ b_{88}^4 + 2 b_{22}^4
                - b_{32}^1 + b_{00}^2 - 2 b_{44}^2 + b_{88}^2\cr
        a_3&= b_{00}^ 3 + b_{77}^3 + b_{66}^3 - 2 b_{22}^3 - b_{88}^3\cr
        a_4&=-b_{41}^1 - b_{32}^1.}$$
This can only be obtained following lengthy calculations involving several 
choices in solving the linear equations generated by taking the exterior
derivative of equations (3.3).

We will continue working with this new Darboux frame
$\{A_0;\tilde A_1,\ldots,\tilde A_4;A_5,\ldots,A_8\}$.  To simplify
the notation and since we no longer have any need for
our original frame, we will denote this new frame by
$\{A_0; A_1,\ldots, A_4;A_5,\ldots,A_8\}$.  Following this change of frame
the matrix $\Omega$ will have the form

$$ \pmatrix{ 
\omega_{00} & \omega_{1} & \omega_{2} & \omega_{3}
 & \omega_{4} & 0 & 0 & 0 & 0 \cr
\omega_{10} & \omega_{11} & \omega_{12} & \omega_{13}
 & \omega_{14} & {{\omega_{1}}\over{2}} & {{\omega_{3}}\over{2}} & 
{{\omega_{4}}\over{2}} & 0 \cr 
\omega_{20} & \omega_{21} & \omega_{22} & -\omega_{14}
 & \omega_{13} & {{\omega_{2}}\over{2}} & {{\omega_{4}}\over{2}} & 
-{{\omega_{3}}\over{2}} & 0 \cr 
\omega_{30} & \omega_{31} & \omega_{32} & \omega_{33}
 & \omega_{34} & 0 & {{\omega_{1}}\over{2}} & -{{\omega_{2}}\over{2}} 
& {{\omega_{3}}\over{2}} \cr 
\omega_{40} & -\omega_{32} & \omega_{31} & \omega_{43}
 & \omega_{44} & 0 & {{\omega_{2}}\over{2}} & {{\omega_{1}}\over{2}} 
& {{\omega_{4}}\over{2}} \cr 
\omega_{50} & \omega_{51} & \omega_{52} & \omega_{53}
 & \omega_{54} & \omega_{55} & \omega_{56} & \omega_{57} 
& 0 \cr 
\omega_{60} & \omega_{61} & \omega_{62} & \omega_{63}
 & \omega_{64} & \omega_{65} & \omega_{66} & \omega_{67} 
& \omega_{68} \cr 
\omega_{70} & \omega_{71} & \omega_{72} & \omega_{73}
 & \omega_{74} & \omega_{75} & \omega_{76} & \omega_{77} 
& \omega_{78} \cr 
\omega_{80} & \omega_{81} & \omega_{82} & \omega_{83}
 & \omega_{84} & 0 & \omega_{86} & \omega_{87} 
& \omega_{88} \cr
}\eqno(3.14) $$
We note that in this matrix we have also included the fact from (3.10)
that $\omega_{85}=\omega_{58}=0$.

We again note that from the linear equations (3.10) and (3.11),
and from taking the exterior derivative of the equations in (3.3) we do
have other information concerning relations between 
the forms in matrix (3.14) but, at this
point we are not able to determine any other forms in terms of the 
basis forms.  To be more precise, only some of the information contained in 
the equations (3.10) and (3.11) is displayed in the equation (3.14). There are
other relations between the forms, but we do not know that any forms can
be written in terms of others.

We also notice that we would like to have the anti-symmetric structure in 
the rectangle with corners (6,6) and (7,7).  Similar to (3.12)
it is possible to make a change of frame of the form
$$\tilde A_j=A_j+\sum_{k=1}^{4}l_{jk} A_k.\eqno(3.15)$$
It is not possible to make a change of frame only involving $A_6$ and $A_7$
that respects Kuiper's normal form.  Therefore we perform the change
of frame given in (3.15) for $j=5,6,7,8$.  We then look for values of $l_{jk}$
which will result in $\ot_{66}=\ot_{77}$, 
as well as preserve (3.10) and (3.11).
A solution for the $l_{jk}$ that satisfy this requirement
exist, and we will henceforth assume that we are 
working with this adjusted frame.

Once we have made these pair of changes of frame, it is possible to 
return to equations (3.3) and their
exterior derivatives, taking into account these new frames.
From this we have the following immediate consequences:
the anti-symmetric structure exists in minors cornered at
(1,1) and (2,2) as well as (3,3) and (4,4).  Further, it is possible
to determine that $\o_{67}=-\o_{76}$.

$$ \pmatrix{ 
\omega_{0,0} & \omega_{1} & \omega_{2} & \omega_{3}
 & \omega_{4} & 0 & 0 & 0 & 0 \cr
\omega_{10} & \omega_{11} & \omega_{12} & \omega_{13}
 & \omega_{14} & {{\omega_{1}}\over{2}} & {{\omega_{3}}\over{2}} & 
{{\omega_{4}}\over{2}} & 0 \cr 
\omega_{20} & -\omega_{12} & \omega_{11} & -\omega_{14}
 & \omega_{13} & {{\omega_{2}}\over{2}} & {{\omega_{4}}\over{2}} & 
-{{\omega_{3}}\over{2}} & 0 \cr 
\omega_{30} & \omega_{31} & \omega_{32} & \omega_{33}
 & \omega_{34} & 0 & {{\omega_{1}}\over{2}} & -{{\omega_{2}}\over{2}} 
& {{\omega_{3}}\over{2}} \cr 
\omega_{40} & -\omega_{32} & \omega_{31} & -\omega_{34}
 & \omega_{33} & 0 & {{\omega_{2}}\over{2}} & {{\omega_{1}}\over{2}} 
& {{\omega_{4}}\over{2}} \cr 
\omega_{50} & \omega_{51} & \omega_{52} & \omega_{53}
 & \omega_{54} & \omega_{55} & \omega_{56} & \omega_{57} 
& 0 \cr 
\omega_{60} & \omega_{61} & \omega_{62} & \omega_{63}
 & \omega_{64} & \omega_{65} & \omega_{66} & \omega_{67} 
& \omega_{68} \cr 
\omega_{70} & \omega_{71} & \omega_{72} & \omega_{73}
 & \omega_{74} & \omega_{75} & -\omega_{67} & \omega_{66} 
& \omega_{78} \cr 
\omega_{80} & \omega_{81} & \omega_{82} & \omega_{83}
 & \omega_{84} & 0 & \omega_{86} & \omega_{87} 
& \omega_{88} \cr
}$$

\bigskip
\penalty-100
\ni{\sl \S 3.7.  Third Change of Frame.}
\penalty200
\medskip

By equations (3.10) we know that 
$\omega_{58}=\omega_{85}=0$.  We can take the 
exterior derivatives of these equations to obtain
$$\eqalign{0&=\o_{81}\w\o_{1}+\o_{82}\w\o_{2}+2\o_{86}\w\o_{65}
        +2\o_{87}\w\o_{75}\cr
        0&=\o_{85}\w\o_3+\o_{54}\w\o_{4}+\o_{54}\w\o_{4}+2\o_{56}\w\o_{68}
                +2\o_{57}\w\o_{78}.}\eqno(3.16)$$

By expanding (3.16) in terms of the basis forms using (3.4), and 
collecting the terms
which are coefficients of wedges of basis forms $\o_l\wedge \o_k$ we arrive
at a system of equations in $b_{jk}^l$.   

Now, using the equations (3.3) and (3.16) we can gather additional
information concerning the forms in $\Omega$ without making any
more changes of frame.  The following can be verified:
$$\eqalign{\o_{65} = \o_{31},\ \ \ \ \ \ \ \ \ &\o_{68}=\o_{13}\cr
        \o_{75} = \o_{41}=-\o_{32},\ \ \ \ \ \ \ \ \ &\o_{78} = \o_{14} \cr
        \o_{56} = 2\o_{13},\ \ \ \ \ \ \ \ \ &\o_{57}=2\o_{14}\cr
        \o_{86} = 2\o_{31},\ \ \ \ \ \ \ \ \ &\o_{87}=2\o_{41}=-2\o_{32}.\cr
        }$$
Consequently, the information that we have concerning our
matrix $\Omega$ is the following:

$$ \pmatrix{ 
\omega_{0,0} & \omega_{1} & \omega_{2} & \omega_{3}
 & \omega_{4} & 0 & 0 & 0 & 0 \cr
\omega_{10} & \omega_{11} & \omega_{12} & \omega_{13}
 & \omega_{14} & {{\omega_{1}}\over{2}} & {{\omega_{3}}\over{2}} & 
{{\omega_{4}}\over{2}} & 0 \cr 
\omega_{20} & -\omega_{12} & \omega_{11} & -\omega_{14}
 & \omega_{13} & {{\omega_{2}}\over{2}} & {{\omega_{4}}\over{2}} & 
-{{\omega_{3}}\over{2}} & 0 \cr 
\omega_{30} & \omega_{31} & \omega_{32} & \omega_{33}
 & \omega_{34} & 0 & {{\omega_{1}}\over{2}} & -{{\omega_{2}}\over{2}} 
& {{\omega_{3}}\over{2}} \cr 
\omega_{40} & -\omega_{32} & \omega_{31} & -\omega_{34}
 & \omega_{33} & 0 & {{\omega_{2}}\over{2}} & {{\omega_{1}}\over{2}} 
& {{\omega_{4}}\over{2}} \cr 
\omega_{50} & \omega_{51} & \omega_{52} & \omega_{53}
 & \omega_{54} & \omega_{55} & 2\omega_{13} & 2\omega_{14} 
& 0 \cr 
\omega_{60} & \omega_{61} & \omega_{62} & \omega_{63}
 & \omega_{64} & \omega_{31} & \omega_{66} & \omega_{67} 
& \omega_{13} \cr 
\omega_{70} & \omega_{71} & \omega_{72} & \omega_{73}
 & \omega_{74} & -\omega_{32} & -\omega_{67} & \omega_{66} 
& \omega_{14} \cr 
\omega_{80} & \omega_{81} & \omega_{82} & \omega_{83}
 & \omega_{84} & 0 & 2\omega_{31} & -2\omega_{32} 
& \omega_{88} \cr
}$$

Next, we know that with the current choice of frame
$\o_{14}+\o_{23}=0$ and $\o_{41}+\o_{32}=0$, see (3.13).  We can
take the exterior derivative
of these equations to get
$$\eqalign{0&=d\o_{14}+d\o_{23}\cr
        0&=d\o_{41}+d\o_{32}.}$$
After solving this system of equations we are ready
to make another change of frame.

The next change of frame is necessary to get zeros in the location 
$\omega_{53},\ \omega_{54}$ and $\omega_{81},\ \omega_{82}$. 
Let
$$\tilde A_5=A_5+a_5 A_0$$
$$\tilde A_8=A_8+a_8 A_0.$$
This change of frame will alter the forms only in the top and
bottom rows of the minor cornered at
(5,0) and (8,5).
Once again, if we let $\ot_{5j}$ and $\ot_{8j}$ for $j=1,\ldots 5$
denote the Maurer-Cartan
forms associated with this new frame, we see that by choosing
$a_5=2b_{10}^1-2b_{74}^1$ and $a_8=2b_{40}^4-2b_{62}^4$
we  get
$$\ot_{52}=\ot_{53}=\ot_{81}=\ot_{82}=0.$$
Following this change of frame the  
matrix $\Omega$ now has the appearance

$$ \pmatrix{ 
\omega_{0,0} & \omega_{1} & \omega_{2} & \omega_{3}
 & \omega_{4} & 0 & 0 & 0 & 0 \cr
\omega_{10} & \omega_{11} & \omega_{12} & \omega_{13}
 & \omega_{14} & {{\omega_{1}}\over{2}} & {{\omega_{3}}\over{2}} & 
{{\omega_{4}}\over{2}} & 0 \cr 
\omega_{20} & -\omega_{12} & \omega_{11} & -\omega_{14}
 & \omega_{13} & {{\omega_{2}}\over{2}} & {{\omega_{4}}\over{2}} & 
-{{\omega_{3}}\over{2}} & 0 \cr 
\omega_{30} & \omega_{31} & \omega_{32} & \omega_{33}
 & \omega_{34} & 0 & {{\omega_{1}}\over{2}} & -{{\omega_{2}}\over{2}} 
& {{\omega_{3}}\over{2}} \cr 
\omega_{40} & -\omega_{32} & \omega_{31} & -\omega_{34}
 & \omega_{33} & 0 & {{\omega_{2}}\over{2}} & {{\omega_{1}}\over{2}} 
& {{\omega_{4}}\over{2}} \cr 
\omega_{50} & \omega_{51} & \omega_{52} & 0
 & 0 & \omega_{55} & 2\omega_{13} & 2\omega_{14} 
& 0 \cr 
\omega_{60} & \omega_{61} & \omega_{62} & \omega_{63}
 & \omega_{64} & \omega_{31} & \omega_{66} & \omega_{67} 
& \omega_{13} \cr 
\omega_{70} & \omega_{71} & \omega_{72} & \omega_{73}
 & \omega_{74} & -\omega_{32} & -\omega_{67} & \omega_{66} 
& \omega_{14} \cr 
\omega_{80} & 0 & 0 & \omega_{83}
 & \omega_{84} & 0 & 2\omega_{31} & -2\omega_{32} 
& \omega_{88} \cr
}\eqno(3.17)$$
We will henceforth assume that we are working with this frame.

\bigskip
\penalty-100
\ni{\sl \S 3.8.  Fourth change of frame}
\penalty200
\medskip

At this point it is also possible to 
check that from the equations already
solved, relationships exist between forms that have 
not yet been discussed.
For example, individual forms can be checked to verify that the 
anti-symmetric form discussed previously exists in other blocks of the matrix
as well. For example, we previously determined that $\Omega$ has
the anti-symmetric structure in the minors with corners at
(1,1) and (2,2), and also at (3,3) and (4,4), and 
finally at (6,6) and (7,7).  
Taking the exterior derivative of this collection of 
equations results in 
$$\eqalign{0&=d\o_{11}-d\o_{22}\cr
          0&=d\o_{12}+d\o_{21}\cr
        0&=d\o_{33}-d\o_{44}\cr
        0&=d\o_{34}+d\o_{43}\cr
        0&=d\o_{66}-d\o_{77}\cr
    0&=d\o_{67}+d\o_{76}.}$$
If these are expanded using (3.4), and we once again collect coefficients
of wedges of basis forms, we obtain a system of linear equations in
terms of the $b_{jk}^l$.

The major region of the matrix that we have not yet examined are the minor 
with corners at $(6,1)$ and $(7,2)$ and the minor with corners $(6,3)$ and
$(7,4)$.  While the latter of these has the anti-symmetric structure
discussed previously, the former does not. 

\medskip
\ni{\bf Remark.}  In the case $k=2$ that we are 
discussing here, the minor with
corners $(6,1)$ and $(7,2)$ appears to have an ``inverted'' anti-symmetric
structure.  In fact, this does not follow through in higher dimensions.  
Examining the matrix for the $k=4$ case will display this
clearly.  See Appendix A for the Maurer-Cartan matrix in the case $k=4$.
\medskip

We will now make a change of frame to get the minor with corners
$(6,3)$ and $(7,4)$  into the form of the matrix (3.1).
The change
of frame that achieves this is
$$\eqalign{\tilde A_6&=A_6+a_6 A_0\cr
           \tilde A_7&=A_7+a_7 A_0.}$$
A change of frame of this type will only alter forms in 
the minor with corners
at (6,0) and (7,4).
We can solve for appropriate choices for
$a_6$ and $a_7$, and see that
$$\eqalign{a_6&=2b_{40}^2-b_{84}^2\cr
        a_7&=2b_{40}^1-b_{84}^1.}$$

Following this change of frame, the matrix $\Omega$ has the appearance

$$ \pmatrix{ 
\omega_{0,0} & \omega_{1} & \omega_{2} & \omega_{3}
 & \omega_{4} & 0 & 0 & 0 & 0 \cr
\omega_{10} & \omega_{11} & \omega_{12} & \omega_{13}
 & \omega_{14} & {{\omega_{1}}\over{2}} & {{\omega_{3}}\over{2}} & 
{{\omega_{4}}\over{2}} & 0 \cr 
\omega_{20} & -\omega_{12} & \omega_{11} & -\omega_{14}
 & \omega_{13} & {{\omega_{2}}\over{2}} & {{\omega_{4}}\over{2}} & 
-{{\omega_{3}}\over{2}} & 0 \cr 
\omega_{30} & \omega_{31} & \omega_{32} & \omega_{33}
 & \omega_{34} & 0 & {{\omega_{1}}\over{2}} & -{{\omega_{2}}\over{2}} 
& {{\omega_{3}}\over{2}} \cr 
\omega_{40} & -\omega_{32} & \omega_{31} & -\omega_{34}
 & \omega_{33} & 0 & {{\omega_{2}}\over{2}} & {{\omega_{1}}\over{2}} 
& {{\omega_{4}}\over{2}} \cr 
\omega_{50} & \omega_{51} & \omega_{52} & 0
 & 0 & \omega_{55} & 2\omega_{13} & 2\omega_{14} 
& 0 \cr 
\omega_{60} & \omega_{61} & \omega_{62} & \omega_{10}
 & \omega_{20} & \omega_{31} & \omega_{66} & \omega_{67} 
& \omega_{13} \cr 
\omega_{70} & \omega_{71} & \omega_{72} & -\omega_{20}
 & \omega_{10} & -\omega_{32} & -\omega_{67} & \omega_{66} 
& \omega_{14} \cr 
\omega_{80} & 0 & 0 & \omega_{83}
 & \omega_{84} & 0 & 2\omega_{31} & -2\omega_{32} 
& \omega_{88} \cr
}$$

\bigskip
\penalty-100
\ni{\sl \S 3.9.  Determining the matrix $\Omega$}
\penalty200
\medskip

At this point, all of the necessary changes of frame have been made, and 
what remains is to check that all the 
information that we have gathered so far is sufficient to 
completely characterize $\Omega$.
The changes of frame that we have performed have given us
several new relations between forms in the matrix, and if we take the
exterior derivative of these equations, we will obtain relations
concerning other forms.  In particular we have the following relations:
$$\eqalign{0&=\o_{13}-\o_{24}=\o_{23}+\o_{14}  \cr
0&= \o_{31}-\o_{42}= \o_{32}+\o_{41}  \cr
0&= \o_{11}-\o_{22}=\o_{12}+\o_{21}     \cr
0&=\o_{33}-\o_{44}= \o_{34}+\o_{43}       \cr
0&=\o_{66}-\o_{77}=\o_{67}+\o_{76}       \cr
0&=\o_{31}-\o_{65}=\o_{41}-\o_{66}      \cr
0&=2\o_{31}-\o_{86}=2\o_{41}-\o_{87}   \cr
0&=\o_{13}-\o_{68}=\o_{14}-\o_{78}     \cr
0&=2\o_{13}-\o_{56}=2\o_{14}-\o_{57}   \cr
0&=\o_{53}=\o_{54}     \cr
0&=\o_{81}=\o_{82}     \cr
0&=\o_{63}-\o_{74}=\o_{73}+\o_{64}     \cr
0&=\o_{61}+\o_{72}=\o_{71}+\o_{62}. }$$

By taking the exterior derivative of this set of equations and
expanding this set of equations in terms of the basis forms (3.4),
and equating the coefficients of wedges of basis forms $\o_j$ to zero, we
get a collection of approximately 150
equations involving the $b_{jk}^l$.  
Solving the equations in the order given makes the calculations convenient
in solving the entire set.
Once this entire set is solved, one
can then examine $\Omega$ entry by entry 
and verify the the matrix $\Omega$ has the same
form as (3.1).

$$ \pmatrix{ 
\omega_{00} & \omega_{1} & \omega_{2} & \omega_{3} & \omega_{4} & 0 & 0
& 0 & 0 \cr 
\omega_{10} & \omega_{11} &  \omega_{12} & \omega_{13} & \omega_{14} &
{\omega_{1}\over 2} & { \omega_{3}\over 2} & {\omega_{4}\over 2} & 0 \cr
\omega_{20} &  -\omega_{12} & \omega_{11} & -\omega_{14} & \omega_{13} &
{\omega_{ 2}\over 2} & {\omega_{4}\over 2} & -{\omega_{3}\over 2} & 0
\cr 
\omega_{30} & \omega_{31} & \omega_{32} & \omega_{33} & \omega_{34}  & 0
& {\omega_{1}\over 2} & -{\omega_{2}\over 2} & {\omega_{3 }\over 2} \cr
\omega_{40} & -\omega_{32} & \omega_{31} & -\omega_{34}  & \omega_{33} &
0 & {\omega_{2}\over 2} & {\omega_{1}\over 2} & {\omega_{4}\over 2} \cr 
0 & 2\,\omega_{10} & 2\,\omega_{20} & 0 & 0 & \omega_{55} &
2\,\omega_{13} & 2\,\omega_{14} & 0 \cr 
0 &  \omega_{30} & \omega_{40} & \omega_{10} & \omega_{20} & \omega_{31}
& \omega_{66} & \omega_{67} & \omega_{13} \cr
0 & \omega_{40} & -\omega_{30} & -\omega_{20} & \omega_{10} &
-\omega_{32} & - \omega_{67} & \omega_{66} & \omega_{14} \cr
0 & 0 & 0 & 2\,\omega_{30} & 2\,\omega_{40} & 0 & 2\,\omega_{31} &
-2\,\omega_{32} &  \omega_{88} \cr
} $$

It now follows from the Lemma in Section 1 that a neighborhood $U$ of
the point $p$ is projectively equivalent to a piece of a standard embedding
of $\C P^2$ in $\P^8$.

\bigskip
\penalty-100
\ni{\sl \S 3.10. Completion of the proof}
\penalty200
\medskip

To finish the proof one has to show that the whole submanifold $M$ is 
projectively equivalent to the standard embedding of $\C P^2$. Let $W$
be the set of interior points of the
connected component containing $p$ of the subset of $M$ on  
which the tight immersion of $M$ agrees with the standard embedding of 
$\C P^2$ after $f$ has been composed with a projective transformation
that brings a neighborhood of $p$ into the standard embedding. We want to
show that $W=M$. Assume the boundary of $W$ is nonempty  and let
$q$ be a boundary point. By continuity, the second fundamental form of $f$
coincides with the one of the standard embedding in $q$. Let $\xi$ be
a normal vector of the standard embedding at $q$ such that the 
corresponding height function of the standard embedding 
has a nondegenerate minimum in $q$. Then
the corresponding height function of $f$ has a nondegenerate minimum
in $q$ since its Hessian at $q$ is ${\rm Hess}(X,Y)=\langle
II(X,Y),\xi\rangle$ and hence the same as the Hessian of the corresponding
height function of the standard embedding at $q$.

Hence there is a neighborhood of $q$ on which $f$
coincides with a standard immersion up to a projective transformation. It 
follows that $f$ is analytic around $q$. Hence $q$ is in the interior of $W$, 
a contradiction. It follows that
$W=M$ and we have finished the proof.

\bigskip
\penalty-100
\ni{\sl \S 3.11.  Some remarks about $k=4$ and $k=8$}
\penalty200
\medskip

For the most part, the method used in the case $k=2$ also works for
the cases $k=4$ and $k=8$, and in fact it can also be used in the case
$k=1$ in place of the method of [GH].  We first compute the matrix of
Maurer-Cartan forms for the standard embedding of $\K P^2$ in
$\P^{3k+2}$, as was done in (3.1).  Next we proceed as in \S3.2 to set
up a matrix of Maurer-Cartan forms for a tight substantial immersion
of highly connected manifold $M^{2k}$ in $\P^{3k+2}$.  In all of the
cases we make use of Kuiper's normal form to determine the upper right
hand block in our of Maurer-Cartan forms as in \S3.2.  The primary
difficulty as for larger values of $k$ is the magnitude and complexity
of the calculations.  In the case when $k=8$ we have a $27\times 27$
matrix $\Omega$ of Maurer-Cartan forms, and the normal form allows us to
describe the 160 forms in the rectangle in the upper right hand corner
with corners $(1,17)$ and $(16,26)$. This is analogous to the
calculations resulting in (3.3).  Proceeding as in \S3.3 we can then take
the exterior derivative of these equations.  Collecting the
coefficients of wedges of the basis forms as in (3.6) results in
approximately 19200 linear equations.  The calculations
in all cases then proceed similarly. There is
effectively no difference between the algorithms used in the
different cases other than limits on sizes of certain arrays, the
sizes of matrices and the complexity of the calculations.

To give some indication of the
degree of similarity between the cases, we have
included the final desired matrix in the quaternion case in
the Appendix.  Once again, the normal form described by Kuiper completely
describes the upper right hand corner of the matrix.  Also, the
anti-symmetric structure that occurs in the $k=2$ case extends
quite predictably to the $k=4$ case.   The only region of the matrix where
anything unpredictable occurs in the case $k=4$ is the minor cornered
at (10,1) and (13,4).  Fortunately, we determine the frame
by determining the forms in the minor (10,5) and (13,8) as
in (3.17) and the
structure of the minor cornered at (10,1) and (13,4) is a consequence.

\vfil\eject
\centerline{\msc Appendix A:}  
\centerline{The matrix of Maurer-Cartan forms for the standard}
\centerline{embedding of $\H P^2$ in $\R P^{14}$.}

\vskip 1truecm
\hskip 4truecm\hbox{\setbox0 \hbox{
\localmag{800}{8}{10}
$\displaystyle \pmatrix{ 
\omega_{0,0} & \omega_{0,1} & \omega_{0,2} & \omega_{0,3}
 & \omega_{0,4} & \omega_{0,5} & \omega_{0,6} & \omega_{0,7} & \omega
_{0,8} & 0 & 0 & 0 & 0 & 0 & 0 \cr 
\omega_{1,0} & \omega_{1,1} & 
\omega_{1,2} & \omega_{1,3} & \omega_{1,4} & \omega_{1,5} & \omega_{1,
6} & \omega_{1,7} & \omega_{1,8} & {\omega_{0,1}\over 2} & {\omega_{0,
5}\over 2} & {\omega_{0,6}\over 2} & {\omega_{0,7}\over 2} & {\omega_{
0,8}\over 2} & 0 \cr 
\omega_{2,0} & -\omega_{1,2} & \omega_{1,1} & 
\omega_{2,3} & \omega_{2,4} & -\omega_{1,6} & \omega_{1,5} & -\omega_{
1,8} & \omega_{1,7} & {\omega_{0,2}\over 2} & {\omega_{0,6}\over 2} & 
-{\omega_{0,5}\over 2} & {\omega_{0,8}\over 2} & -{\omega_{0,7}\over 2
} & 0 \cr 
\omega_{3,0} & -\omega_{1,3} & -\omega_{2,3} & \omega_{1,1}
 & \omega_{3,4} & -\omega_{1,7} & \omega_{1,8} & \omega_{1,5} & -
\omega_{1,6} & {\omega_{0,3}\over 2} & {\omega_{0,7}\over 2} & -{
\omega_{0,8}\over 2} & -{\omega_{0,5}\over 2} & {\omega_{0,6}\over 2}
 & 0 \cr 
\omega_{4,0} & -\omega_{1,4} & -\omega_{2,4} & -\omega_{3,4}
 & \omega_{1,1} & -\omega_{1,8} & -\omega_{1,7} & \omega_{1,6} & 
\omega_{1,5} & {\omega_{0,4}\over 2} & {\omega_{0,8}\over 2} & {\omega
_{0,7}\over 2} & -{\omega_{0,6}\over 2} & -{\omega_{0,5}\over 2} & 0
 \cr 
\omega_{5,0} & \omega_{5,1} & \omega_{5,2} & \omega_{5,3} & 
\omega_{5,4} & \omega_{5,5} & \omega_{5,6} & \omega_{5,7} & \omega_{5,
8} & 0 & {\omega_{0,1}\over 2} & -{\omega_{0,2}\over 2} & -{\omega_{0,
3}\over 2} & -{\omega_{0,4}\over 2} & {\omega_{0,5}\over 2} \cr 
\omega
_{6,0} & -\omega_{5,2} & \omega_{5,1} & -\omega_{5,4} & \omega_{5,3}
 & -\omega_{5,6} & \omega_{5,5} & \omega_{6,7} & \omega_{6,8} & 0 & {
\omega_{0,2}\over 2} & {\omega_{0,1}\over 2} & -{\omega_{0,4}\over 2}
 & {\omega_{0,3}\over 2} & {\omega_{0,6}\over 2} \cr 
\omega_{7,0} & -
\omega_{5,3} & \omega_{5,4} & \omega_{5,1} & -\omega_{5,2} & -\omega_{
5,7} & -\omega_{6,7} & \omega_{5,5} & \omega_{7,8} & 0 & {\omega_{0,3}
\over 2} & {\omega_{0,4}\over 2} & {\omega_{0,1}\over 2} & -{\omega_{0
,2}\over 2} & {\omega_{0,7}\over 2} \cr 
\omega_{8,0} & -\omega_{5,4}
 & -\omega_{5,3} & \omega_{5,2} & \omega_{5,1} & -\omega_{5,8} & -
\omega_{6,8} & -\omega_{7,8} & \omega_{5,5} & 0 & {\omega_{0,4}\over 2
} & -{\omega_{0,3}\over 2} & {\omega_{0,2}\over 2} & {\omega_{0,1}
\over 2} & {\omega_{0,8}\over 2} \cr 
0 & 2\,\omega_{1,0} & 2\,\omega_{
2,0} & 2\,\omega_{3,0} & 2\,\omega_{4,0} & 0 & 0 & 0 & 0 & \omega_{9,9
} & 2\,\omega_{1,5} & 2\,\omega_{1,6} & 2\,\omega_{1,7} & 2\,\omega_{1
,8} & 0 \cr 
0 & \omega_{5,0} & \omega_{6,0} & \omega_{7,0} & \omega_{8
,0} & \omega_{1,0} & \omega_{2,0} & \omega_{3,0} & \omega_{4,0} & 
\omega_{5,1} & \omega_{10,10} & \omega_{10,11} & \omega_{10,12} & 
\omega_{10,13} & \omega_{1,5} \cr 
0 & \omega_{6,0} & -\omega_{5,0} & -
\omega_{8,0} & \omega_{7,0} & -\omega_{2,0} & \omega_{1,0} & \omega_{4
,0} & -\omega_{3,0} & -\omega_{5,2} & -\omega_{10,11} & \omega_{10,10}
 & \omega_{11,12} & \omega_{11,13} & \omega_{1,6} \cr 
0 & \omega_{7,0}
 & \omega_{8,0} & -\omega_{5,0} & -\omega_{6,0} & -\omega_{3,0} & -
\omega_{4,0} & \omega_{1,0} & \omega_{2,0} & -\omega_{5,3} & -\omega_{
10,12} & -\omega_{11,12} & \omega_{10,10} & \omega_{12,13} & \omega_{1
,7} \cr 
0 & \omega_{8,0} & -\omega_{7,0} & \omega_{6,0} & -\omega_{5,0
} & -\omega_{4,0} & \omega_{3,0} & -\omega_{2,0} & \omega_{1,0} & -
\omega_{5,4} & -\omega_{10,13} & -\omega_{11,13} & -\omega_{12,13} & 
\omega_{10,10} & \omega_{1,8} \cr 
0 & 0 & 0 & 0 & 0 & 2\,\omega_{5,0}
 & 2\,\omega_{6,0} & 2\,\omega_{7,0} & 2\,\omega_{8,0} & 0 & 2\,\omega
_{5,1} & -2\,\omega_{5,2} & -2\,\omega_{5,3} & -2\,\omega_{5,4} & 
\omega_{14,14} \cr}
$
}\rotl0}

\vfil\eject

\ni{\bf References}

\font\rrm=cmr8
\font\rsl=cmsl8
\font\rbf=cmbx8
\rrm

\parskip 5pt
\baselineskip 9pt plus 2pt minus 2pt
\bigskip

\item{[BS]} R. Bott and H. Samelson, {\rsl Applications of the theory of 
Morse to symmetric spaces}, Amer. J. Math. {\rbf 80} (1958) 964--1029.

\item{[CL1]} S.-S. Chern and R. Lashof, {\rsl On the total curvature of
immersed manifolds}, Amer. J. Math. {\rbf 79} (1957), 306--318. 

\item{[CL2]} S.-S. Chern and R. Lashof, {\rsl On the total curvature of
immersed manifolds, II}, Michigan Math. J. {\rbf 5} (1958), 5--12.

\item{[CR]} T. E. Cecil and P. J. Ryan, {\rsl Tight and taut immersions
of manifolds}, Research Notes in Math., vol. 107, Pitman Publ., Boston, MA, 1985

\item{[GH]} Ph. Griffiths and J. Harris, {\rsl Algebraic geometry and
local differential geometry}, Ann. Scient. \'Ec. Norm. Sup., IV. Ser. 
{\rbf 12} (1979), 355--452.

\item{[H]} A. Hurwitz, {\rsl \"Uber die Komposition der quadratischen
Formen},
Math. Ann. {\rbf 88} (1923), 1--25.  Also in 
{\rsl A. Hurwitz: Mathematische Werke II}, 
pp. 641--666, Verlag von Emil Birkh\"auser \& Cie., Basel, 1933.

\item{[KT]} S. Kobayashi and M. Takeuchi, {\rsl Minimal imbeddings of 
R-spaces}, J. Differential Geom. {\rbf 2} (1968), 203--215.

\item{[K1]} N. H. Kuiper, {\rsl Sur les immersions \`a courbure totale
minimale}, S\'eminaire de Topologie et de G\'eom\'etrie Differentielle
C. Ehresmann, vol.~II, Paris, 1961.

\item{[K2]} N. H. Kuiper, {\rsl On convex maps}, Nieuw Archief voor Wisk. 
{\rbf 10} (1962), 147--164.

\item{[K3]} N. H. Kuiper, {\rsl Minimal total absolute curvature 
for immersions}, Invent. Math. {\rbf 10} (1970), 209--238.

\item{[K4]} N. H. Kuiper, {\rsl Tight embeddings and maps.  Submanifolds
of geometrical class three in ${\rbf E}^n$}, The Chern Symposium 1979 (Proc. 
Internat. Sympos., Berkeley, Calif., 1979) pp. 97--145, Springer-Verlag,
Berlin, Heidelberg, New York, 1980. 



\item{[LP]} J. A. Little and W. F. Pohl, {\rsl On tight immersions of 
maximal codimension}, Invent. Math. {\rbf 13} (1971), 179--205.

\item{[Sa]} T. Sasaki, {\rsl On the Veronese embedding and related system
of differential equations}, 
Global differential geometry and global analysis (Proc. Conf. Berlin, 1990),
pp. 210--247, Lecture Notes in Math. 1481, Springer-Verlag, 
Berlin, Heidelberg, New York, 1991. 

\item{[S]} F. Severi, {\rsl Intorno ai punti doppi impropri di una
superficie generale dello spazio a quattro dimensioni, 
e a suoi punti tripli apparenti}, Rendiconti di Palermo {\rbf 15} (1901), 
33--51.  Also in {\rsl Memorie Scelte} (Selected Works), Vol. I, pp. 19--38,
Dott. Cesare Zuffi -- Editore, Bologna, 1950.


\item{[T]} S.-S. Tai, {\rsl Minimum imbeddings of compact 
symmetric spaces of rank
one}, J. Differential Geom. {\rbf 2} (1968), 55--66. 

\item{[Th]} G. Thorbergsson, {\rsl Tight immersions of highly connected 
manifolds}, Comment. Math. Helv. {\rbf 61} (1986), 102--121.


\vskip 2truecm

\font\scaps=cmcsc10 at 8pt
\hbox{\parindent=0pt\parskip=0pt
\vbox{\hsize=3.7truein
\obeylines
{\scaps
Ross Niebergall
Mathematics and Computer Science
University of Northern British Columbia
Prince George, British Columbia, Canada}
{\sl rossn@unbc.edu}
}  
\hskip 1.5truecm
\vbox{\hsize=2.7truein
\obeylines
{\scaps
Gudlaugur Thorbergsson
Mathematisches Institut
Universit\"at zu K\"oln
K\"oln, Germany}
{\sl gthorbergsson@mi.uni-koeln.de}
}
}

\end